\documentclass[11pt]{elsarticle}

\usepackage{graphicx}
\usepackage{amssymb,xfrac}
\usepackage{amsthm,amsmath,dsfont}
\usepackage{thmtools}

\usepackage{lineno}

\usepackage[colorlinks=true, citecolor=blue]{hyperref}
\usepackage[algoruled,titlenotnumbered]{algorithm2e}
\usepackage{algorithmic}
\algsetup{linenosize=\small}

\usepackage{wrapfig}
\usepackage{natbib}
\usepackage{comment}
\usepackage{multirow}
\usepackage{amsmath,bm}
\usepackage{graphicx}
\usepackage[lofdepth,lotdepth]{subfig}

\usepackage{float}
\usepackage{subfig}
\usepackage{wrapfig}
\usepackage{framed,xcolor}
\usepackage{newfloat}
\usepackage[xcolor,framemethod=TikZ]{mdframed}
\usepackage{stmaryrd}
\usepackage{mathtools}
\usepackage{mathpazo}

\usepackage[colorinlistoftodos]{todonotes}

\usepackage{longtable}

\newcommand{\RR}{\mathbb{R}}

\newcommand{\dd}{\mathrm{d}}
\newcommand{\vot}{\mathrm{vot}}

\newtheoremstyle{theoremdd}
{\topsep}
{\topsep}
{\itshape}
{0pt}
{\fontfamily{cmss}\selectfont\bfseries}
{.}
{ }
{\thmname{#1}\thmnumber{ #2}\thmnote{ (#3)}}

\theoremstyle{theoremdd}

\theoremstyle{definition}

\newproof{pf}{Proof}

\mdfdefinestyle{myStyle}{roundcorner=5pt,backgroundcolor=brown!20,linecolor=red!40!black,linewidth=1.5pt}
\mdfdefinestyle{myStyle1}{roundcorner=2pt,innertopmargin=10pt,linewidth=1.0pt}

\graphicspath{{./Figures/}}

\biboptions{authoryear}
\usepackage{titlesec}
\titleformat*{\section}{\fontfamily{cmss}\selectfont\large\bfseries}
\titleformat*{\subsection}{\fontfamily{cmss}\selectfont\normalsize\bfseries}
\titleformat*{\subsubsection}{\fontfamily{cmss}\selectfont\normalsize}

\begin{document}	
\begin{frontmatter}

	\title{\fontfamily{cmss}\selectfont Pay to change lanes: A cooperative lane-changing strategy for connected/automated driving}
	
	\author[add1]{Dianchao Lin}
	\author[add1]{Li Li}
	\author[add1,add2]{Saif Eddin Jabari\corref{cor1}}
	\cortext[cor1]{Corresponding author, e-mail: \url{sej7@nyu.edu}
	}
	
	\address[add1]{Tandon School of Engineering, New York University, Brooklyn, New York, U.S.A.}
	\address[add2]{Division of Engineering, New York University Abu Dhabi, Saadiyat Island, P.O. Box 129188, Abu Dhabi, U.A.E.}
	
	{ \fontfamily{cmss}\selectfont\large\bfseries		
		\begin{abstract}
			{ \normalfont\normalsize
This paper proposes a cooperative lane changing strategy using a transferable utility games framework. This allows vehicles to engage in transactions where gaps in traffic are created in exchange for monetary compensation.  The proposed approach is best suited to discretionary lane change maneuvers.  We formulate gains in travel time, referred to as time differences, that result from achieving higher speeds. These time differences, coupled with value of time, are used to formulate a utility function, where utility is transferable.  We also allow for games between connected vehicles that do not involve transfer of utility.  We apply Nash bargaining theory to solve the latter.  A cellular automaton is developed and utilized to perform simulation experiments that explore the impact of such transactions on traffic conditions (travel-time savings, resulting speed-density relations and shock wave formation) and the benefit to vehicles. The results show that lane changing with transferable utility between drivers can help achieve win-win results, improve both individual and social benefits without resulting in any adverse effects on traffic characteristics in general and, in fact, result in slight improvement at traffic densities outside of free-flow and (bumper-to-bumper) jammed traffic.
			}
		\end{abstract}
	}
	
	
\end{frontmatter} 

\section{Introduction}
\label{S:intro}

Lane-changing is one of the fundamental maneuvers in vehicular traffic dynamics.  Vehicles change lanes to achieve desired speeds (discretionary lane-changing), to avoid unsafe conditions or to move into turning/exit lanes (mandatory lane changing).  A majority of models of both types of lane change maneuvers describe them as discrete decision processes carried out by vehicles that are considering/attempting to change lanes \citep{gipps1986model,kesting2007general,zheng2014recent,kamal2015efficient,du2015autonomous,keyvan2016categorization,li2016lane,pan2016modeling,bevly2016lane}.    We refer to \citep{ahmed1996models,toledo2003modeling} for a classical reference on discrete choice methods for lane changing and to \citep{pan2016modeling} for a more recent background on these types of lane changing models. These tend to ignore the competition for space that may arise between vehicles and how this competition affects their decisions.  This has given rise to game theoretic techniques in modeling lane-changing dynamics \citep{kita1999merging,kita2002game,wang2015game,meng2016dynamic,liu2007game,talebpour2015modeling,oyler2016game,li2017game}. A typical setting in these approaches is one in which a discrete set of maneuvers (typically two or three) are being considered by a vehicle that is attempting to change lanes (the \textit{target vehicle}) and a vehicle that is in the target lane but behind the target vehicle (the \textit{lag vehicle}).  For example, the target vehicle may have the choice set \{change lane, do not change lane\} while the lag vehicle has the choice set \{give way, do not give way\} \citep{rahman2013review}. 

The different approaches in the literature vary in how they model the payoffs associated with pure strategies, which may vary depending on whether the maneuver is mandatory or discretionary. Some papers only consider lane changing games for mandatory behavior, such as merging \citep{kita1993effects, kita1999merging, pei2006control}. Most game-theoretic approaches consider lane changing to be non-cooperative games, the outcomes of which are either Nash or Stackelberg equilibria depending on how the game is modeled \citep{yoo2013stackelberg,li2016hierarchical,yu2018human}. Cooperative strategies have also been considered recently \citep{wang2015game,yao2017optimizing,zimmermann2018carrot}. A common feature of the latter is that vehicles are assumed to be selfless; one in which \textit{cooperative} vehicles (typically under some form of control) will take actions that maximize the collective or group utility, not their own.  This leads to winners \textit{and} losers.

Automation and vehicle to vehicle (V2V) communication present an opportunity to re-think lane-changing strategies.  These allow vehicles to broadcast their payoffs, which can vary from vehicle to vehicle and for the same vehicle from trip to trip, i.e., depending on trip purpose \citep{hossan2016investigating}.  Communication also allows vehicles to engage in bargaining (and/or repeated) games.  These two features culminate in a departure from the traditional (non-cooperative) game-theoretic lane-changing approaches in which decisions are made without communication\footnote{Although it is common to assume that the vehicles are perfectly knowledgeable (of the payoffs).}.  Indeed, with vehicle to vehicle (V2V) communication capabilities, connected vehicles can easily engage in transactions based on their individual travel needs. For example, quick mobile payment without transaction costs has gained popularity in China.

In light of this, we propose a discretionary lane changing paradigm, suitable for an automated world, in which gaps in traffic are envisaged as resources (or goods) that can be traded.  In simple terms, we propose a lane changing mechanism that allows vehicles to purchase right of way or compensate other vehicles for allowing them to change lanes.  From a modeling stand-point, we propose modeling lane changing as transferable utility (TU) games with side payments \citep{thomas2008game, myerson2013game}.  Our approach also allows for vehicles to refuse to engage in TU games as well; in this case, we consider Nash bargaining.  It can be shown (when there are no transaction fees) that the outcomes of these games are at least pareto efficient \citep{coase2013problem}.  To the best of our knowledge, this is the first time TU games are applied to lane changing dynamics.   

This paper is organized as follows: \autoref{S:model} describes the problem setting and formulates the utility functions and the game's payoffs in \autoref{SS:classification} - \autoref{SS:time_difference}.  The remainder of \autoref{S:model} presents the lane change games with, transfer of utility, side payments (\autoref{SS:TUgame}), and games between connected vehicles that do not wish to engage in transactions, i.e., games with non-transferable utility (\autoref{SS:nonTUgame}). \autoref{S:simulation} presents a numerical example and a set of simulation experiments to test the proposed model, analyze the results from the aspects of cost-effectiveness and impact on traffic flow. \autoref{S:Conc} concludes the paper.

\section{Methodology}
\label{S:model}

\subsection{Problem description}
\label{SS:classification}
We consider discretionary lane change maneuvers.  The game setting we assume is one that is played between pairs of vehicles, a lane changing (target) vehicle and a lag vehicle in the target lane. We shall focus on this simple setting, but the ideas are generalizable to games involving more than two vehicles or settings involving two lane changing vehicles.  For example, the proposed approach can be easily extended to situations with multiple pairwise lane change games: one simply parallelizes the proposed approach.  A more general framework is one where a vehicle engaged in transactions with more than one vehicle simultaneously, e.g., a vehicle that wishes to make two immediate lane change maneuvers.  The first maneuver can be analyzed in a similar way to the pairwise approach described in this paper.  Pricing the second maneuver is more subtle as it depends on the outcome of the first maneuver.  However, the proposed approach can still be utilized as a building block for such sophisticated scenarios.

\begin{figure}[h!]
	\centering
	\resizebox{0.80\textwidth}{!}{%
		\includegraphics{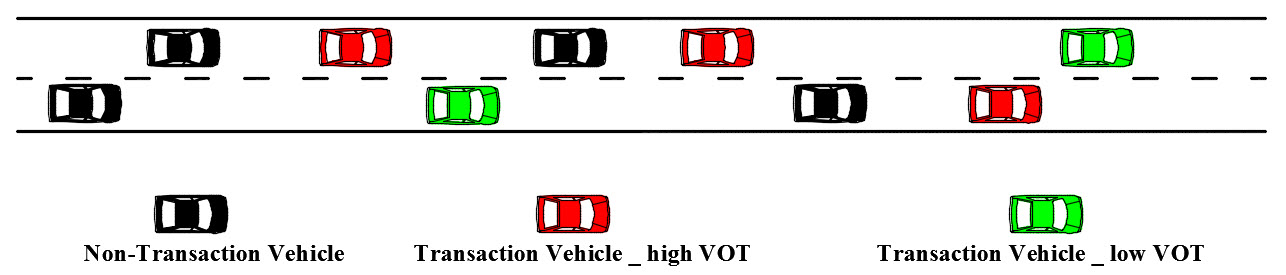}}
	\caption{Transaction vehicles (TVs) and non-transaction vehicles (NTVs).} 
	\label{F:different_vehicles}
\end{figure}

It is assumed that vehicles involved in the game can communicate positions, speeds, accelerations, and values of time (VOTs).  This paper assumes that vehicles communicate these variables \textit{truthfully} and leave the interesting question of untruthfulness to future research.  We assume mixed traffic: vehicles capable of (and willing to) engage in \textit{lane change transactions}, referred to as \textit{transaction vehicles} (TV), and those that either do not possess transaction capabilities (or do not wish to engage in transactions), which are referred to as \textit{non-transaction vehicles} (NTV); see \autoref{F:different_vehicles}.  We assume that vehicles equipped with communication capabilities are free to choose between being TVs or NTVs.  We also assume some level of automation in the TVs allowing for minimal input from drivers over the course of a trip, in which multiple such lane change decisions (transactions) may take place.

The types of games considered in this paper are those in which utility can be \textit{transfered} (or traded) between vehicles.  While the motivations for changing lanes can include numerous factors, such as comfort, safety, and speed, not all of these factors can be traded.  For example, safety cannot (and should not) be traded.  Utility can only be transfered when a common currency that is valued equally by both vehicles is used.  It is for this reason that we formulate the utility function below using VOT. 
The mathematical notation used in this paper is listed in \autoref{t1}.
\newpage
{ \small
\begin{longtable}{|c p{0.85\linewidth}|} 
\caption{Notation} \label{t1}\\
\hline
Variable & Description \\
\hline
	$v_1^{i}$ & the higher speed choice (of vehicle $i$) \\
	$v_2^{i}$ & the lower speed choice (of vehicle $i$) \\
	$v_{\mathrm{E}}^{i}$ &  the expected equilibrium speed (of vehicle $i$) \\
	$v^i$ & the speed choice of vehicle $i$ \\
	$a_1^{i}$ & the acceleration from $v_1$ to $v_{\mathrm{E}}$ (of vehicle $i$) \\
	$a_2^{i}$ & the acceleration from $v_2$ to $v_{\mathrm{E}}$ (of vehicle $i$) \\
	$a_{\mathrm{pos}}$ & the value of $a_1$ and $a_2$ we used in simulation when they are positive \\
	$a_{\mathrm{neg}}$ & the value of $a_1$ and $a_2$ we used in simulation when they are negative \\
	$t_a$ & the average time it takes for a vehicle to complete a lane-change \\
	$t_{b1}$ & the time when the speed of vehicle changes from $v_1$ to $v_{\mathrm{E}}$ \\
	$t_{b2}$ & the time when the speed of vehicle changes from $v_2$ to $v_{\mathrm{E}}$ \\
	$t_b$ & $\max\{t_{b1}, t_{b2}\}$ \\
	$S_a$ & the difference in distance achieved when choosing $v_1$ over $v_2$ from time 0 to $t_a$ \\
	$S_b$ &  the difference in distance achieved when choosing $v_1$ over $v_2$ from time $t_a$ to $t_b$ \\
	$S$ & the difference in distance achieved when choosing $v_1$ over $v_2$ from time 0 to $t_b$, $S=S_a + S_b$ \\
	$t_{\dd}^{i}$ & the time difference between achieving $v_1$ and $v_2$ (for vehicle $i$) \\
	$c_{\vot}^{i}$ & the coefficient representing the VOT (of vehicle $i$) \\
	$u^i$ & the utility of vehicle $i$ \\
	$\mathbf{A}$, $\mathbf{B}$ & the utility matrices of vehicles A and B, respectively \\
	$p$ & the probability for vehicle A to change lanes in its threat strategy in a TU game \\
	$q$ & the probability for vehicle B to not give way in its threat strategy in a TU game \\
	$T_{\mathrm{A}}, T_{\mathrm{B}}$ & the payoffs of vehicles A and B, respectively, at their threat strategies \\
	$Q_{\mathrm{A}}, Q_{\mathrm{B}}$ & the payoffs of vehicles A and B, respectively, at status quo \\
	$\omega^*$ & the total maximal achievable utility by vehicles A and B in a TU game \\
	$(i^*,j^*)$ & the strategy pair that achieves the maximal payoff $\omega^*$, also the final decision of TU game \\
	$\sigma$ & the side payment in a TU game \\
	$N_{\mathrm{A}},N_{\mathrm{B}}$ & the payoff of vehicle A and B at the Nash bargaining solution in an NTU game \\
	$(i',j')$ & the final decision of NTU game \\
	$p_{\mathrm{sd}}$ & the probability for vehicle to slow down in the simulation\\
	$v_{\max}$ & the maximum speed of vehicle in the simulation \\
	$\beta$ & the benefit index of vehicles \\
	\hline
\end{longtable}
}

\subsection{Utility function}
\label{SS:utilityFunction}
\begin{figure}[h!]
	\centering
	\resizebox{0.4\textwidth}{!}{%
		\includegraphics{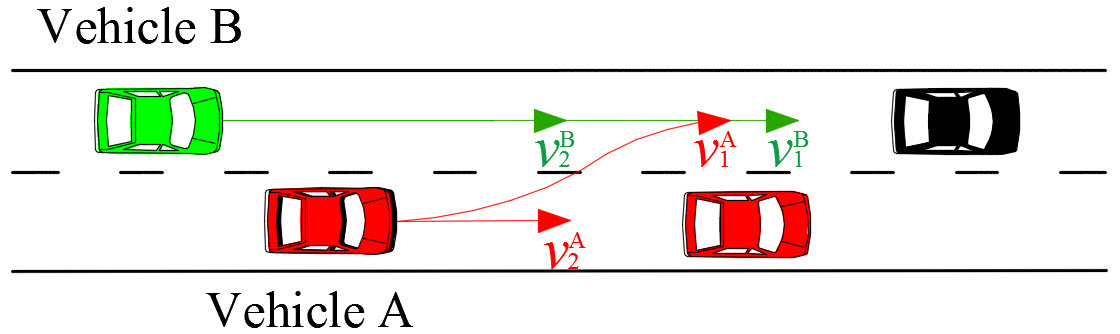}}
	\caption{Lane change scenarios between two TVs.} 
	\label{F:lane_change}
\end{figure}

Consider two vehicles, A and B, in the lane-changing game depicted in \autoref{F:lane_change}. Vehicle A has the choice set: \{change lanes, do not change lanes\} and vehicle B's choice set is \{give way, do not give way\}.    
The payoffs are summarized in \autoref{t2}; $A_{ij}$ and $B_{ij}$ are the payoffs to vehicles A and B, respectively,  associated with actions $i$ and $j$. 
\begin{table}[h!]
	\centering
	\caption{Utility matrix for lane changing game}
	\begin{tabular}{|c|c|c|c|}
		\hline
		\multicolumn{2}{|c|}{\multirow{2}{*}{Actions}} &
		\multicolumn{2}{|c|}{vehicle B} \\
		\cline{3-4}
		\multicolumn{2}{|c|}{} & Do not give way & Give way \\
		\hline
		\multirow{2}{*}{vehicle A}
		& Change lanes & ($A_{11}$,$B_{11}$) & ($A_{12}$,$B_{12}$) \\
		\cline{2-4}
		& Do not change lanes & ($A_{21}$,$B_{21}$) & ($A_{22}$,$B_{22}$) \\
		\hline
	\end{tabular}
	\label{t2}
\end{table}
If A chooses to change lanes and B chooses \textit{not} to give way, we assume that both vehicles get a large negative payoff (e.g., they collide).  If A changes lanes and B gives way, we denote the speeds achieved by vehicles A and B by $v^{\mathrm{A}}_1$ and $v^{\mathrm{B}}_2$, respectively.  If A stays as a result of B not giving way, we denote the speeds achieved by A and B by $v^{\mathrm{A}}_2$ and $v^{\mathrm{B}}_1$, respectively.  Here, we have that $v^{\mathrm{A}}_1 > v^{\mathrm{A}}_2$ and $v^{\mathrm{B}}_1 > v^{\mathrm{B}}_2$.  When A chooses to stay despite B giving way, they both assume the lower speeds $v^{\mathrm{A}}_2$ and $v^{\mathrm{B}}_2$. 

For vehicle A, the utility function, denoted $u^{\mathrm{A}}$ ($u^{\mathrm{B}}$ for B) is related to both its own speed choice $v^{\mathrm{A}}$, and that of vehicle B denoted $v^{\mathrm{B}}$.  Similarly, $u^{\mathrm{B}}$ is related to both $v^{\mathrm{B}}$ and $v^{\mathrm{A}}$.  Vehicle A's utility function is given by	
 \begin{align}
	u^{\mathrm{A}}(v^{\mathrm{A}},v^{\mathrm{B}}) = \begin{cases} 
	-M & v^{\mathrm{A}} = v^{\mathrm{A}}_1, v^{\mathrm{B}} = v^{\mathrm{B}}_1 \\
	c_{\vot}^{\mathrm{A}} t_{\dd}^{\mathrm{A}} & v^{\mathrm{A}} = v^{\mathrm{A}}_1, v^{\mathrm{B}} = v^{\mathrm{B}}_2\\
	0 & v^{\mathrm{A}} = v^{\mathrm{A}}_2\\
	\end{cases},
	\label{E:utility}
	\end{align} 
where $M$ is a large positive number, $c_{\vot}^{\mathrm{A}}$ is a coefficient capturing VOT of A, and $t_{\dd}^{\mathrm{A}}$ is referred to as \textit{time difference} between choosing lower speed $v^{\mathrm{A}}_2$ and higher speed $v^{\mathrm{A}}_1$ for A.  Using time difference as a means of calculating the utility, the reference point (a.k.a. \textit{datum}) of $u^{\mathrm{A}}$ coincides with the action $v^{\mathrm{A}} = v^{\mathrm{A}}_2$; that is $u^{\mathrm{A}}(v^{\mathrm{A}}_2,\cdot) = 0$.  The latter can be interpreted as: choosing/maintaining the lower speed comes with zero utility.  Vehicle B's utility, $u^{\mathrm{B}}$, can be calculated in a similar fashion. 
The resulting payoff matrices, denoted $\mathbf{A}$ and $\mathbf{B}$ for players A and B, respectively are given as follows
\begin{equation}
\mathbf{A} \equiv \begin{bmatrix} A_{11} & A_{12} \\ A_{21} & A_{22} \end{bmatrix} = \begin{bmatrix} -M & c_{\vot}^{\mathrm{A}} t_{\dd}^{\mathrm{A}} \\0 & 0 \end{bmatrix} \text{ and } \mathbf{B} \equiv \begin{bmatrix} B_{11} & B_{12} \\ B_{21} & B_{22} \end{bmatrix} = \begin{bmatrix} -M & 0 \\ c_{\vot}^{\mathrm{B}} t_{\dd}^{\mathrm{B}} & 0 \end{bmatrix}. \label{E:t3}
\end{equation}
Time difference is the travel time saved with the higher-speed choice over a short time period. How it is calculated is described in \autoref{SS:time_difference} below. As a result of $t_{\dd}$ being defined as \textit{time gained} when compared to the lower speed scenario, the pay-off pair $(A_{22},B_{22}) = (0,0)$.  This is elaborated further below.  It is plausible that the case where B gives way, but A does not change lanes results in vehicle B being annoyed, in which case $B_{22}$ would be negative.  However, we have no way to quantify this loss in utility due to being annoyed and thus ignore it in this paper.

\subsection{Modeling of time difference}
\label{SS:time_difference}
Assume a vehicle is traveling with longitudinal speed $v = v_0$ and at time $t = 0$ an opportunity presents itself for the vehicle to achieve a higher speed (e.g., via a lane change).  We make no assumptions about lane preference in this paper.  That is, vehicles can use any lane to pass.  Consequently, it is reasonable to assume that the two lanes have similar traffic conditions (on average); that is, the (macroscopic) traffic densities are the same and, hence, that the equilibrium speeds do not depend on lane. However, the equilibrium speeds can vary by vehicle depending on their VOT (in addition to traffic density). We do \textit{not} assume that the gaps in traffic are the same and this is what creates the speed gain opportunities. The speed gains occur over a short period of time equal to the length of an equilibration process.  For example, a vehicle changing lanes accelerates to a higher speed then decelerates to the equilibrium speed.  Similarly, a vehicle giving way might decelerate for a short period of time and then accelerate to the equilibrium speed.

\begin{figure}[h!]
	\centering
	\subfloat[][$t$ from 0 to $t_{a}$]{
		\includegraphics[width=0.4\textwidth]{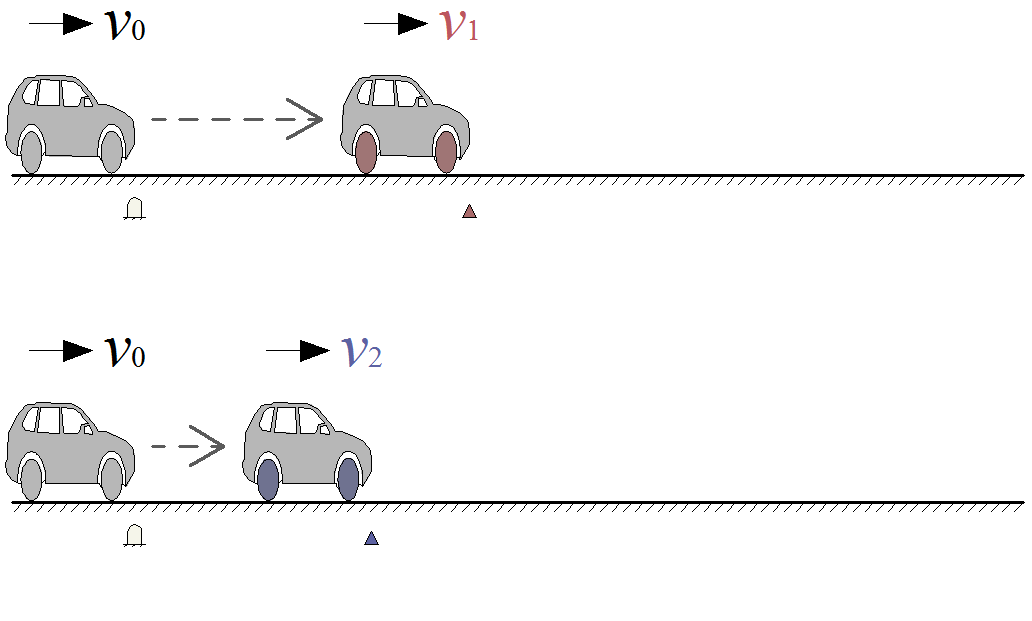}
		\label{F:speed_eventa}} 
	\subfloat[][$t$ from $t_{a}$ to $t_{b}$]{
		\includegraphics[width=0.4\textwidth]{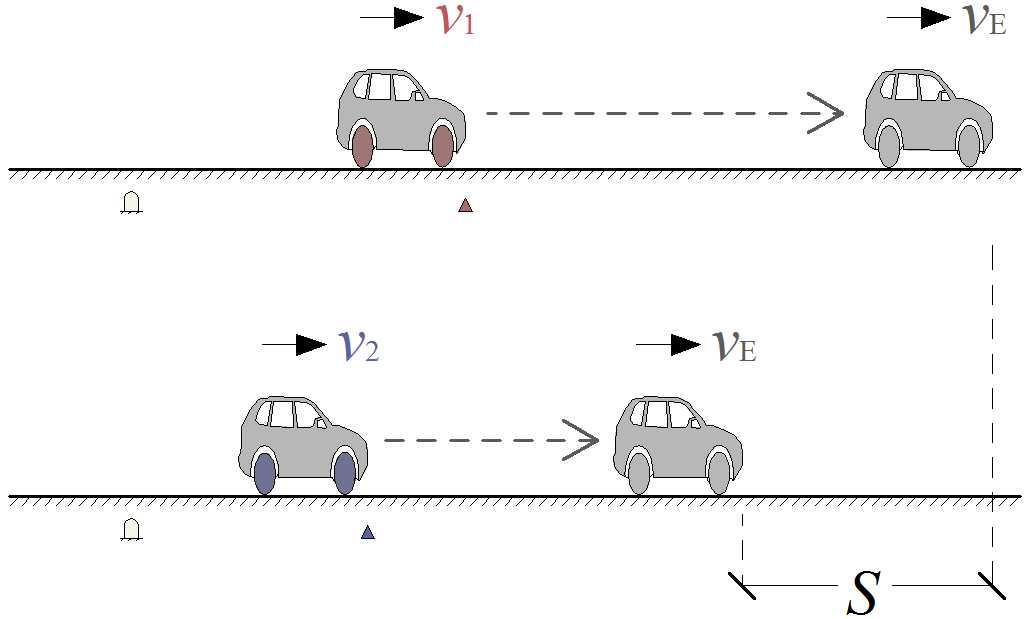}
		\label{F:speed_eventb}}
	\caption{Influence of different speed choices at time $t=0$.}
	\label{F:speed_event}
\end{figure}

\begin{figure}[h!]
	\centering
	\subfloat[Example trajectories for A]{
		\includegraphics[width=0.4\textwidth]{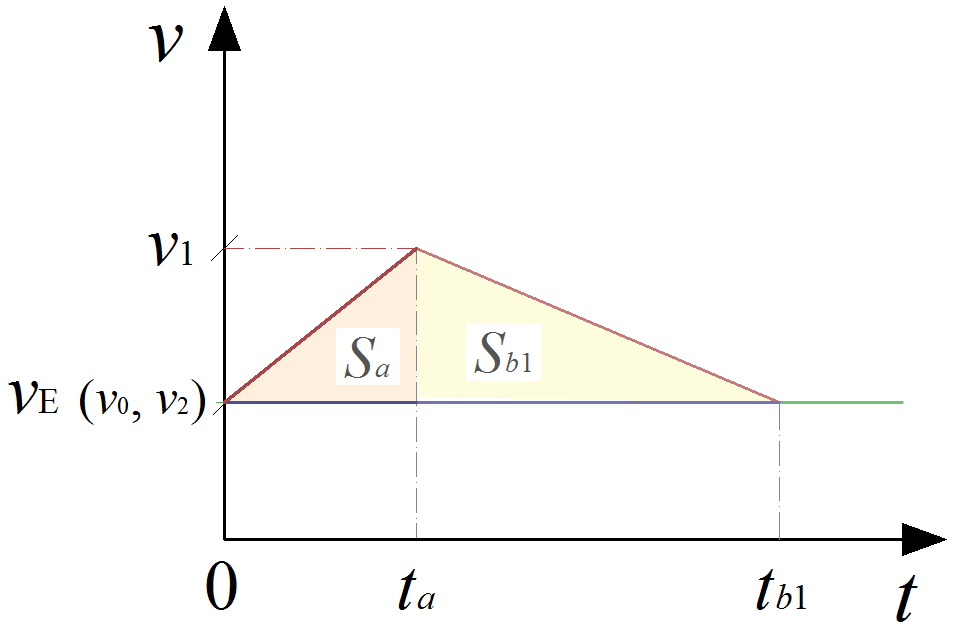}
		\label{F:S1_vehicleA}} 
	\subfloat[Example trajectories for B]{
		\includegraphics[width=0.4\textwidth]{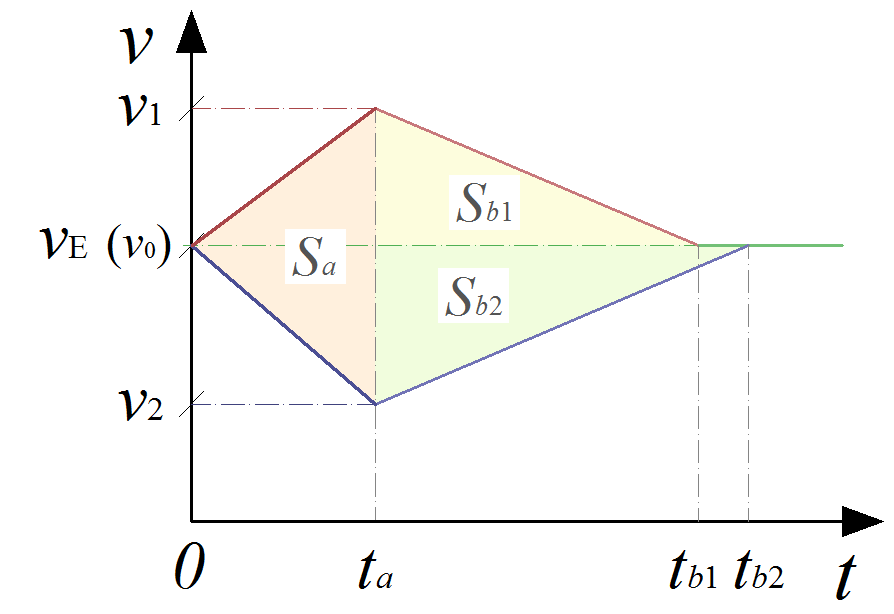}
		\label{F:S1_vehicleB}}
	\caption{Illustration of $S_{a}$, $S_{b1}$ and $S_{b2}$.}
	\label{F:S1}
\end{figure}

Both the target vehicle (vehicle A) and the lag vehicle (vehicle B) can be described as vehicles that are changing their speeds to one of two possible speeds: either $v = v_1$, or $v = v_2$, where $v_1 > v_2$. We denote by $t_a > 0$ the time instant at which the vehicle achieves its new speed.  We assume, without loss of generality, that $t_a$ is the same in both scenarios.  If the higher speed is adopted, illustrated in the top part of \autoref{F:speed_event}\subref{F:speed_eventa}, the vehicle reaches the faster speed $v_1$ over a longer travel distance. If the slower speed is chosen, illustrated in the bottom part of \autoref{F:speed_event}\subref{F:speed_eventa}, it reaches $v_2$ over a shorter travel distance.  Note that $v_0 = v_2$ is allowed but not necessary as depicted in \autoref{F:S1}\subref{F:S1_vehicleA}.  
Assume these two scenarios are identical in all other aspects (including traffic conditions, vehicle performance and characteristics, etc.) so that in both scenarios the vehicle achieves the equilibrium speed eventually.  We denote the equilibrium speed by $v_{\mathrm{E}}$.  
The equilibrium speed $v_{\mathrm{E}}$ is mainly dependent on traffic density. Assume, without loss of generality, that $v_1 \ge v_{\mathrm{E}} \ge v_2$.  Let $t_b$ denote the time instant at which the vehicle achieves the equilibrium speed $v_{\mathrm{E}}$, $t_{b1}$ in the first scenario and $t_{b2}$ in the second scenario. If the vehicle chooses the higher speed, it will have covered a longer distance by time $t_b = \max\{t_{b1},t_{b2}\}$ as shown in \autoref{F:speed_event}\subref{F:speed_eventb}.  We denote the difference in distance covered during the equilibration process by $S$.
If denote the trajectory of the vehicle by $v(\cdot)$, which in the high speed scenario is denoted by $v_{\mathrm{high}}(\cdot)$ and by $v_{\mathrm{low}}(\cdot)$ in the lower speed scenario. Then $S$ is defined as
\begin{equation}
	S \equiv \int_0^{t_b} \big|v(t) - v_{\mathrm{low}}(t)\big| \dd t \label{E:Sdef}
\end{equation}
with $v(t_b) = v_{\mathrm{E}}$. If $v(\cdot) = v_{\mathrm{low}}(\cdot)$, corresponding to the case where vehicle A does not change lanes or when vehicle B gives way, we have that $S = 0$.  Since $t_{\dd} \propto S$ (defined formally below) and $u^{\mathrm{A}}(v_2^{\mathrm{A}},\cdot) \propto t_{\dd}$ (also, $u^{\mathrm{B}}(\cdot,v_2^{\mathrm{B}}) \propto t_{\dd}$), we have that $A_{21} = A_{22} = 0$ and $B_{12} = B_{22} = 0$ as given in \eqref{E:t3}.

In the case where $v(\cdot) = v_{\mathrm{high}}(\cdot)$, $S$ is the \textit{gain in distance} under the higher speed. In this case, under our two speed assumption, \eqref{E:Sdef} can be calculated using simple geometry.  \autoref{F:S1} provides two example trajectory differences, one for vehicle A, \autoref{F:S1}\subref{F:S1_vehicleA}, and one for vehicle B, \autoref{F:S1}\subref{F:S1_vehicleB}.  In both figures, the solid red curve is the vehicle's trajectory in the high-speed scenario, the solid blue curve is their trajectory in the low-speed scenario, and the green line is equilibrium speed.

To calculate $S$, we proceed as follows: we first split $S$ into two parts, $S_a$ and $S_b$, such that $S = S_a + S_b$.  Here, $S_a$ denotes the distance difference from time $0$ to time $t_a$; $S_b$ is the difference in distance covered from time $t_a$ to time $t_b$, which we divide into two parts $S_{b1}$ and $S_{b2}$, such that $S_b = S_{b1} + S_{b2}$.  In \autoref{F:S1}\subref{F:S1_vehicleA}, vehicle A finds a gap in the adjacent lane at time $t = 0$. Their speed at time zero, $v_0$, is assumed to be equal to $v_{\mathrm{E}}$ in the example. If A changes lanes, they temporarily achieve the higher speed of $v_1$.  If $v_1 > v_{\mathrm{E}}$, they then decelerate to $v_{\mathrm{E}}$. If A does not change lanes, they maintain the lower speed of $v_2 = v_{\mathrm{E}}$. We again note that, in general, $v_0$ need not be the equilibrium speed, which would correspond to the case where an equilibration process is interrupted by a lane change maneuver. 
When $v_2 = v_{\mathrm{E}}$ as in \autoref{F:S1}\subref{F:S1_vehicleA}, $S_{b2} = 0$. In \autoref{F:S1}\subref{F:S1_vehicleB}, vehicle B receives a lane-change request from vehicle A. If B gives way, they decelerates to $v_2$ and then accelerates to $v_{\mathrm{E}}$.  If B does not give way, they may temporarily accelerate to $v_1 > v_{\mathrm{E}}$ to close the gap, but must then decelerate to $v_{\mathrm{E}}$.  

From  \autoref{F:S1}\subref{F:S1_vehicleA} and \autoref{F:S1}\subref{F:S1_vehicleB}, we have that $S_a$ is given by
\begin{equation}
	S_a = \frac{(v_1 - v_2) t_a}{2}.
	\label{E:Sa}
\end{equation}
Note that this encompasses the case where vehicle A's speed at time $t=0$ is not the equilibrium speed: $v_0 \ne v_{\mathrm{E}}$.  Similarly, it encompasses the case where vehicle B does not decelerate.  For $S_b$, the times $t_{b1}$ and $t_{b2}$ are generally not known, but it is safe to assume a constant acceleration/deceleration rate that falls in appropriate ranges.  Let $a > 0$ denote this acceleration rate, then
\begin{equation}
	S_{b1} = \frac{(v_{\mathrm{E}} - v_1)^{2}}{-2a}, \label{E:Sb1_1}
\end{equation}
where the negative sign is added to the denominator since $S_{b1}$ describes a deceleration scenario and
\begin{align}
	S_{b2} = \frac{(v_{\mathrm{E}} - v_2)^{2}}{2a} \label{E:Sb2_1}.
\end{align}
Combining \eqref{E:Sa}, \eqref{E:Sb1_1}, and \eqref{E:Sb2_1}, i.e., $S = S_a + S_{b1} + S_{b2}$, we get:
\begin{equation}
	S =  \frac{1}{2} \left[ (v_1 - v_2) t_a + \frac{(v_{\mathrm{E}} - v_1)^{2}}{-a} + \frac{(v_{\mathrm{E}} - v_2)^{2}}{a}\right].
\label{E:S}
\end{equation}

It is worth noting that although \eqref{E:S} was derived under the assumption that $v_1 \ge v_{\mathrm{E}} \ge v_2$, it also applies to other scenarios, $v_1 \ge v_2 \ge v_{\mathrm{E}}$ and $v_{\mathrm{E}} \ge v_1 \ge v_2$. Moreover, the case $v_1 = v_2 = $ the lower speed corresponds to the case where vehicle A does not change lanes (or vehicle B gives way).  Note that, in this case, it follows from \eqref{E:S} that $S=0$.

Finally, we define the time difference $t_{\dd}$ as the average amount of travel time saved during the equilibration process as a result of achieving the higher speed $v_1$ at time $t=0$.  That is:
\begin{equation}
	t_{\dd} \equiv \frac{S}{v_{\mathrm{E}}} =  \frac{1}{2 v_{\mathrm{E}}} \left[ (v_1 - v_2) t_a + \frac{(v_{\mathrm{E}} - v_1)^{2}}{-a} + \frac{(v_{\mathrm{E}} - v_2)^{2}}{a}\right],
\label{E:td}
\end{equation}  
where the parameters $t_a$, and $a$ need to be calibrated.  The speeds $v_1$, $v_2$ and $v_{\mathrm{E}}$ are variables in different TU games, where $v_{\mathrm{E}}$ depends on traffic density in this paper. 
Note that for a road with heterogeneous traffic conditions across lanes, one can generalize the proposed approach to one with different equilibrium speeds in the two lanes and the areas depicted in \autoref{F:S1} would have more components.  We leave this to future research.


\subsection{Transfer of utility and side payments}
\label{SS:TUgame}
In a game with transferable utility (TU), side payments from one vehicle to another are allowed. These are games that involve two (connected) vehicles that wish to engage in a transaction.  This section describes utility transfer and how to calculate the side payment.  Let $\sigma \in \RR$ denote the side payment made by vehicle A to vehicle B.  When $\sigma$ is negative, this is interpreted as a positive payment that is made by vehicle B to vehicle A. When $\sigma < 0$, B does not give way to A, but \textbf{also} makes a payment to A as a result.  This is an important factor distinguishing games where vehicles agree to engage in transactions from those where a transaction does not take place (the non-transaction games described in the next section).

The main idea behind transfer of utility is that through side payment the highest total payoff, denoted $\omega^*$, can be achieved.  Here
\begin{equation}
\omega^* \equiv \underset{i,j \in \{1,2\}\times \{1,2\}}{\max} \big(A_{ij} + B_{ij} \big). \label{E:omega}
\end{equation}
Define the strategy pair that achieves the maximal payoff by
\begin{equation}
(i^*,j^*) \equiv \underset{i,j \in \{1,2\}\times \{1,2\}}{\arg \max} \big(A_{ij} + B_{ij} \big).
\end{equation}
The payoffs achieved this way are $A_{i^*,j^*} - \sigma$ for A and $B_{i^*,j^*} + \sigma$ for B; clearly the total utility is $\omega^*$.  

\begin{figure}[h!]
	\centering
	\subfloat[][Feasible set for the TU game]{
		\includegraphics[width=0.5\textwidth]{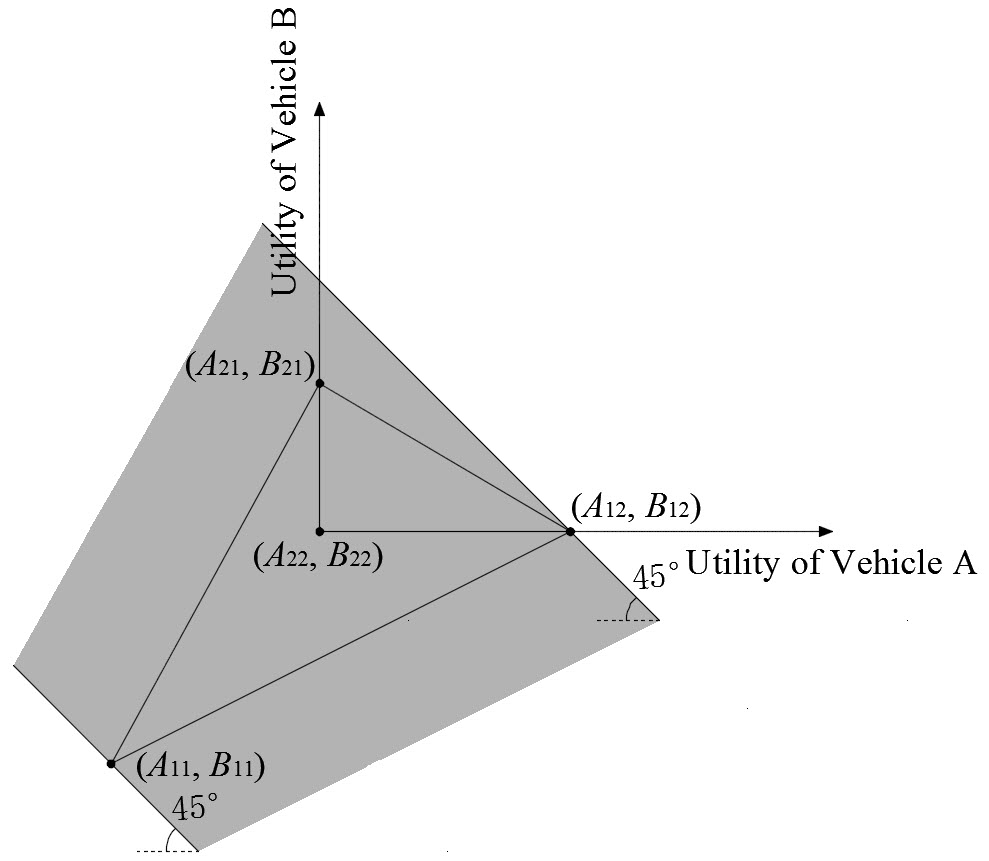}
		\label{F:fs_set_TU}} 
	\subfloat[][Threat point and TU solution]{
		\includegraphics[width=0.5\textwidth]{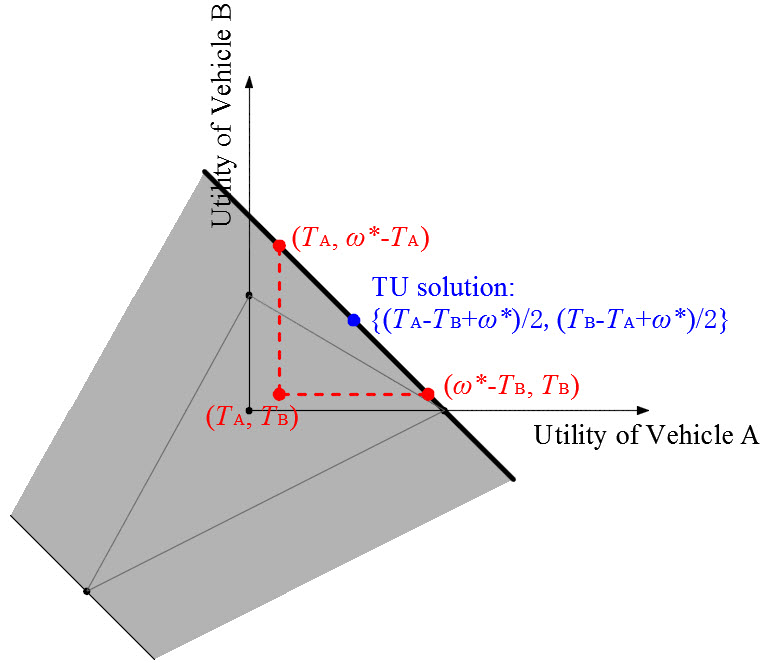}
		\label{F:TU_so}}
	\caption{Feasible set and optimal solution for a TU game.}
	\label{F:TU_game}
\end{figure}

For any $(i,j)$ pair, define $\widetilde{A}_{ij} \equiv A_{ij} - \sigma$ and $\widetilde{B}_{ij} \equiv B_{ij} + \sigma$.  Since, $\widetilde{A}_{ij} + \widetilde{B}_{ij} = \omega_{ij}$, where $\omega_{ij}$ is a constant\footnote{These constants can vary from one $(i,j)$ pair to another and $\omega_{i^*j^*} = \omega^*$.}, we have that the set of payoff pairs associated with the strategy pair $(i,j)$ in a TU game fall along a line of slope -1 that goes through the point $(A_{ij},B_{ij})$.  As such, the set of possible payoffs, the feasible region, for a TU game, denoted $\Omega_{\mathrm{TU}}$, is defined as the convex hull off all possible strategy pairs $(\widetilde{A}_{ij}, \widetilde{B}_{ij})$: 
\begin{equation}
\Omega_{\mathrm{TU}} \equiv \mathsf{Conv} \big( \{(\widetilde{A}_{ij}, \widetilde{B}_{ij}): i,j =1,2 \} \big).
\end{equation}
The feasible region associated with a TU game (the set $\Omega_{\mathrm{TU}}$) is shown in \autoref{F:TU_game}\subref{F:fs_set_TU}.  

In order to determine what the appropriate side payment is, one needs to first investigate the strategy applied in the absence of an agreement.  Such a strategy is known as the \textit{threat strategy}.  Define $\mathbf{p} \equiv [p \quad 1-p]^{\top}$ and $\mathbf{q} \equiv [q \quad 1-q]^{\top}$, where $p$ is the probability that vehicle A chooses to change lanes under their threat strategy and $q$ is the probability that vehicle B chooses to not give way under their threat strategy.  Under their threat strategies, the expected payoffs to vehicles A and B are given, respectively, by $T_{\mathrm{A}} \equiv \mathbf{p}^{\top} \mathbf{A} \mathbf{q}$ and $T_{\mathrm{B}} \equiv \mathbf{p}^{\top} \mathbf{B} \mathbf{q}$. The pair $(T_{\mathrm{A}},T_{\mathrm{B}})$ is known as the \textit{threat point}.

The expected payoffs $T_{\mathrm{A}}$ and $T_{\mathrm{B}}$ (associated with the threat strategies) can be achieved without agreement.  Then, in an agreement the payoffs should be no less than $T_{\mathrm{A}}$ for vehicle A and no less than $T_{\mathrm{B}}$ for vehicle B.  Since the total payoff from the TU solution is $\omega^*$, we have that the TU solution lies between the two points $(T_{\mathrm{A}}, \omega^* - T_{\mathrm{A}})$ and $(\omega^* -T_{\mathrm{B}}, T_{\mathrm{B}})$.  For example, if $(i^*,j^*) = (1,2)$, the range of the TU solution is depicted in  \autoref{F:TU_game}\subref{F:TU_so}.  \cite{thomas2008game} suggests use of the midpoint as a ``natural compromise'': the two vehicles lose equally if the agreement is broken.  The midpoint solution is also depicted in  \autoref{F:TU_game}\subref{F:TU_so}.  

The payoffs associated with the midpoint solution are given by the pair $\big(\frac{T_{\mathrm{A}} - T_{\mathrm{B}} + \omega^*}{2} , \frac{T_{\mathrm{B}} - T_{\mathrm{A}} + \omega^*}{2} \big)$ and we immediately see that A's threat strategy should be chosen in a way that maximizes $T_{\mathrm{A}} - T_{\mathrm{B}}$, while B's should be chosen in a way that minimizes it.  Since
\begin{equation}
T_{\mathrm{A}} - T_{\mathrm{B}} = \mathbf{p}^{\top} \mathbf{A} \mathbf{q} - \mathbf{p}^{\top} \mathbf{B} \mathbf{q} = \mathbf{p}^{\top} (\mathbf{A} - \mathbf{B}) \mathbf{q},
\end{equation}
we have that selecting the threat strategy can be described as a zero-sum game, where A's expected payoff is $\mathbf{p}^{\top} (\mathbf{A} - \mathbf{B}) \mathbf{q}$ and B's expected payoff is $\mathbf{p}^{\top} (\mathbf{B} - \mathbf{A}) \mathbf{q}$. From \eqref{E:t3}, the matrix $\mathbf{A} - \mathbf{B}$ has the following structure
\begin{align}
	\mathbf{A} - \mathbf{B} = \begin{bmatrix}
		0 & c_{\vot}^{\mathrm{A}} t_{\dd}^{\mathrm{A}} \\ -c_{\vot}^{\mathrm{B}} t_{\dd}^{\mathrm{B}} & 0
	\end{bmatrix}.
\end{align}
Since $c_{\vot}^{\mathrm{A}} t_{\dd}^{\mathrm{A}} \ge 0$ and $c_{\vot}^{\mathrm{B}} t_{\dd}^{\mathrm{B}} \ge 0$, we have that $\max_i \min_j (A_{i,j} - B_{i,j}) = \min_j \max_i (A_{i,j} - B_{i,j}) = 0$ and we always have the saddle point \{change lanes, do no give way\}.  That is, the solution is a pure strategy, which corresponds to vehicle A choosing to change lanes with probability $p=1$ and vehicle B choosing to not give way with probability $q = 1$.  While this strategy results in a crash, a crash will not take place since this is a TU game: one in which a total utility of $\omega^*$ is guaranteed.  We hence have that $\mathbf{p} = \mathbf{q} = [1 \quad 0]^{\top}$ and consequently $T_{\mathrm{A}} = T_{\mathrm{B}} = 0$.

\textbf{Side payment}. The payoffs associated with a TU solution (representing the natural compromise): $\frac{1}{2}(T_{\mathrm{A}} - T_{\mathrm{B}} + \omega^*)$ for A and $\frac{1}{2}(T_{\mathrm{B}} - T_{\mathrm{A}} + \omega^*)$ for B.  The side payment is immediately given by 
\begin{align}
\sigma = A_{i^*j^*} - \frac{T_{\mathrm{A}} - T_{\mathrm{B}} + \omega^*}{2} = A_{i^*j^*} - \frac{\omega^*}{2} \label{E:sidepayment}
\end{align}
If $\sigma > 0$, the side payment is made from vehicle A to vehicle B and A changes lanes. If $\sigma < 0$, the side payment is made from vehicle B to vehicle A and B does not give way.  Note that, since $A_{i^*j^*} = \omega^* - B_{i^*j^*}$, we have from \eqref{E:sidepayment} that
\begin{align}
\sigma = \omega^* - B_{i^*j^*} - \frac{T_{\mathrm{A}} - T_{\mathrm{B}} + \omega^*}{2} = - \Big( B_{i^*j^*} - \frac{\omega^*}{2} \Big). \label{E:sidepayment1}
\end{align}
Hence, had we elected to calculate the side payment based on vehicle B's payoffs, we obtain the same side payment with a negative sign, signifying that, in this case, the payment is made in the opposite direction.

\subsection{Games with non-transferable utility}
\label{SS:nonTUgame}
When vehicle A encounters a vehicle that does not wish to engage in a transaction, side payments are not possible.  However, when the two vehicles can communicate, bargaining is possible.  The motivation for considering such situations is that we wish to allow for scenarios in which vehicles do not wish to make side payments and those where it is acknowledged that payment is not always guaranteed for those that wish to receive side payments.  We note that, from a methodological perspective, this does not preclude scenarios that do not involve utility transfer.

\begin{figure}[h!]
	\centering
	\subfloat[][Feasible set for a NTU game]{
		\includegraphics[width=0.5\textwidth]{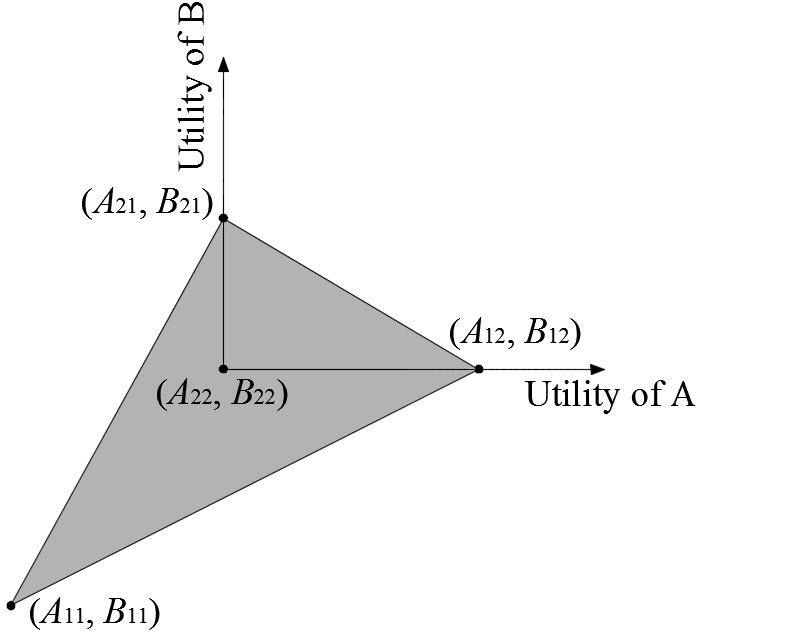}
		\label{F:fs_set_NTU}} 
	\subfloat[][Threat point and Nash bargaining solution]{
		\includegraphics[width=0.5\textwidth]{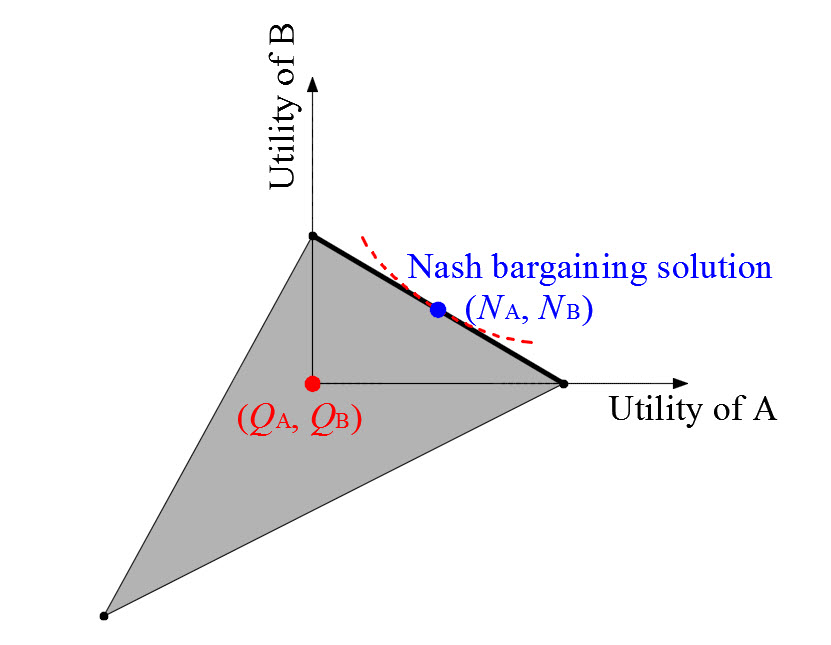}
		\label{F:NTU_so}}
	\caption{Feasible set and optimal solution for a non-transaction game.}
	\label{F:NTU_game}
\end{figure}

The feasible set for a non-transaction lane-change game, denoted by $\Omega_{\mathrm{NTU}}$, is the convex hull of the 4 points, $(A_{ij},B_{ij})_{i,j \in \{1,2\} \times \{1,2\}}$.  That is 
\begin{equation}
	\Omega_{\mathrm{NTU}} \equiv \mathsf{Conv} \big( \{(A_{ij}, B_{ij}): i,j =1,2 \} \big).
\end{equation}
The feasible region is illustrated in \autoref{F:NTU_game}\subref{F:fs_set_NTU}. Communication between vehicle motivates a Nash bargaining solution \citep{nash1950bargaining,nash1953two} for the NTU game. 
The Nash bargaining solution, which we denote by $(N_{\mathrm{A}},N_{\mathrm{B}})$, is the unique solution to the maximization problem
\begin{equation}
	\max \big\{(u_{\mathrm{A}} - Q_{\mathrm{A}})(u_{\mathrm{B}} - Q_{\mathrm{B}}) :  (u_{\mathrm{A}},u_{\mathrm{B}}) \in \Omega_{\mathrm{NTU}}, u_{\mathrm{A}} \ge Q_{\mathrm{A}}, u_{\mathrm{B}} \ge Q_{\mathrm{B}} \big\}, \label{E:NTU}
\end{equation}
where $(Q_{\mathrm{A}},Q_{\mathrm{B}})$ is known as the \textit{status-quo point}. The status quo point occurs before an agreement; $Q_{\mathrm{A}}$ and $Q_{\mathrm{B}}$ are the utilities achieved by vehicles A and B, respectively, if they do not play the game.  It is hence natural to set $(Q_{\mathrm{A}},Q_{\mathrm{B}}) = (0,0)$ corresponding to the strategy \{A does not change lanes, B gives way\}.  Note that in the literature, $(Q_{\mathrm{A}},Q_{\mathrm{B}})$ is sometimes referred to as a threat point; we avoid the latter nomenclature to avoid confusion between $(Q_{\mathrm{A}},Q_{\mathrm{B}})$ and $(T_{\mathrm{A}},T_{\mathrm{B}})$ in the TU game above.

Returning to the bargaining problem \eqref{E:NTU}, we have that
\begin{equation}
(N_{\mathrm{A}},N_{\mathrm{B}}) = \arg \max \big\{ u_{\mathrm{A}} u_{\mathrm{B}} : (u_{\mathrm{A}},u_{\mathrm{B}}) \in \Omega_{\mathrm{NTU}}, u_{\mathrm{A}} \ge 0, u_{\mathrm{B}} \ge 0 \big\}. \label{E:NTU1}
\end{equation}
The contour lines of the objective function in \eqref{E:NTU1} are curves that increase in value the farther one moves away from the origin $(0,0)$.  Since the solution lies in the positive quadrant, the optimal solution lies along the line connecting the two points $(A_{21},B_{21})$ and $(A_{12},B_{12})$, depicted in \autoref{F:NTU_game}\subref{F:NTU_so} as a solid black line.  The equation of this line is $u_{\mathrm{B}} = (-B_{21}/A_{12}) u_{\mathrm{A}} + B_{21}$.  Substituting this into the objective function, we have that 
\begin{align}
	N_{\mathrm{A}} = \frac{1}{2} A_{12} ~ \mbox{ and } ~ N_{\mathrm{B}} = \frac{1}{2} B_{21}.
	\label{E:pA}
\end{align}
Denote the final decision in NTU game as $(i', j')$. We interpret \eqref{E:pA} as: the outcome of the game is either $(i', j') = (1, 2)$, meaning \{A changes lanes, B gives way\} \textbf{or} $(i', j') = (2, 1)$, meaning \{A does not change lanes, B does not give way\}, each outcome with probability 0.5.

\subsection{Model summary and simulation}
\begin{figure}[h!]
	\centering
	\resizebox{0.8\textwidth}{!}{%
		\includegraphics{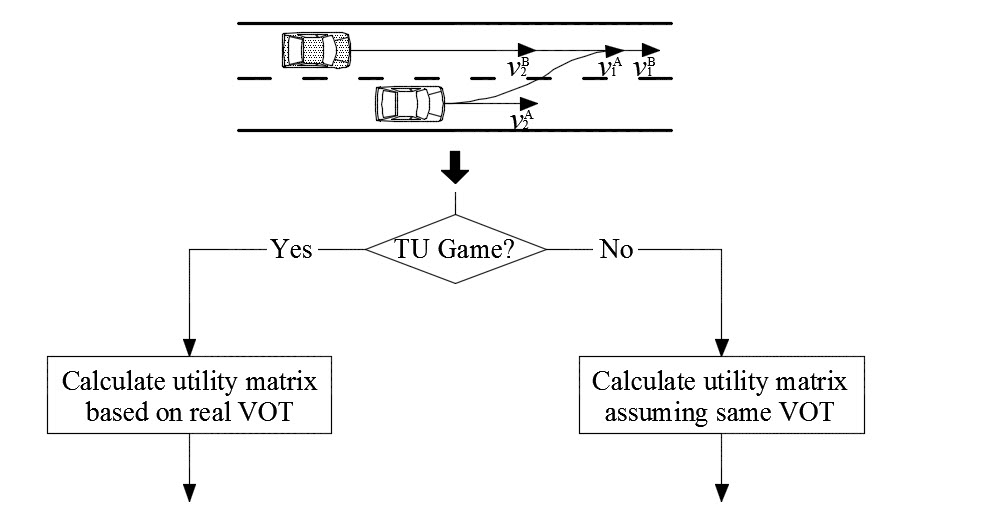}}
	\resizebox{0.8\textwidth}{!}{%
		\includegraphics{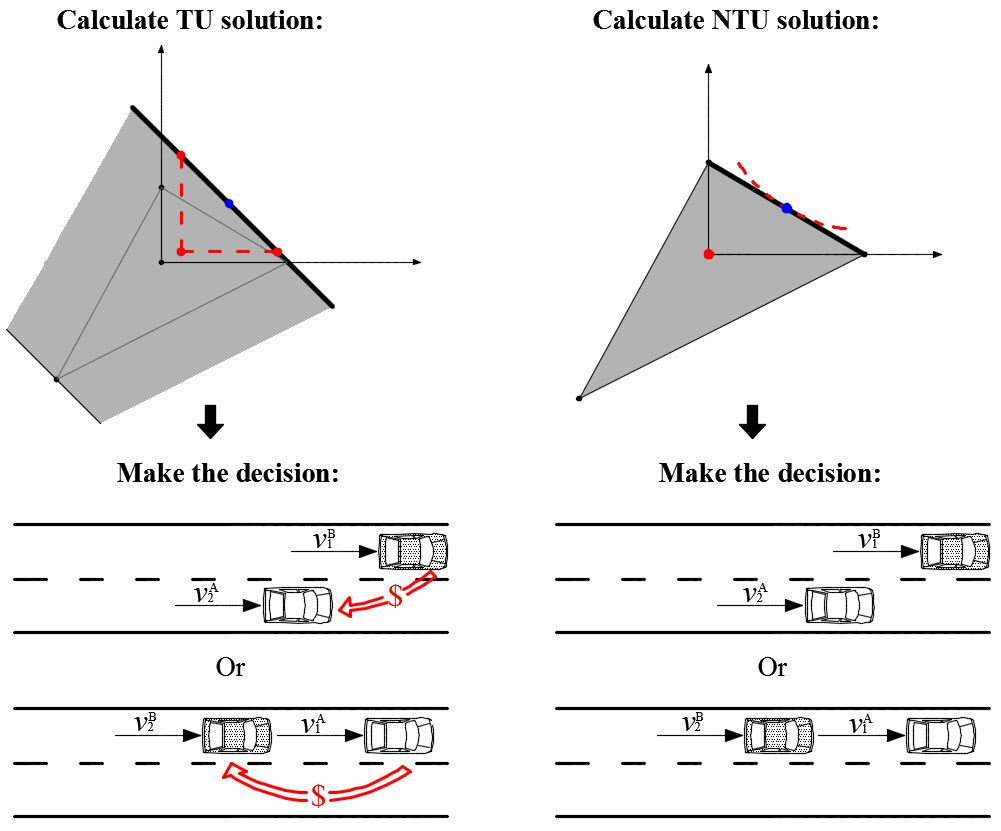}}
	\caption{Lane change decision process flow chart.} 
	\label{F:flowchart}
\end{figure}

The overall process is summarized in \autoref{F:flowchart}. 
In situations where vehicles are not completely aware of the utilities of other vehicles, one can extend the present approach summarized in \autoref{F:flowchart} to include games such as those in \citep{talebpour2015modeling} (assuming no cooperation). 
Safety considerations can also be appended to the present framework, where mandatory lane change conditions arise: see for example the models presented in \citep[Chapter 14]{treiber2013traffic}.  To simulate and test the lane changing game, we propose a cellular automaton \citep{nagel1992cellular, maerivoet2005cellular} given in \autoref{A:1} below.
{
	\small
	\begin{algorithm}[H]
		\caption{Simulation of Traffic Dynamics}
		\label{A:1}
		\begin{algorithmic}
			\STATE \textbf{Input}:
			\STATE \Indp Cell size, number of cells, discrete time step length, inflow rate, 
			\STATE probability of slow-down ($p_{\mathrm{sd}}$), $\{c_{\vot}\}$, percentage of different vehicle classes, max speed ($v_{\max}$), total simulation time, initial speeds and positions of vehicles.
			\STATE \Indm \textbf{Iterate}:
			\STATE \textbf{For} each time step \textbf{do}:
			\STATE \Indp \textbf{For} each vehicle, starting from downstream most \textbf{do}:
			\STATE \Indp \underline{Acceleration}: if $v < v_{\max}$, $v \mapsfrom v + 1$.
			\STATE $d_{\mathrm{s}} \mapsfrom $ distance to leader in subject lane (\# cells)  
			\STATE $d_{\mathrm{t}} \mapsfrom $ distance to leader in target lane (\# cells)  
			\STATE \textbf{If} strategy chosen is stay \textbf{then} 
			\STATE \Indp $v_{\mathrm{stay}} \mapsfrom \min(v, \left \lceil \frac{d_{\mathrm{s}}-1}{2} \right \rceil)$.
			\STATE \Indm \textbf{End If}
			\STATE \textbf{If} strategy chosen is change lanes \textbf{then} 
			\STATE \Indp $v_{\mathrm{change}} \mapsfrom \min(v, v_{\mathrm{stay}} + 1, \left \lceil \frac{d_{\mathrm{t}}-1}{2} \right \rceil)$.
			\STATE \Indm \textbf{End If}
			\STATE \underline{Randomization}: With probability $p_{\mathrm{sd}}$, $v_{\mathrm{stay}} \mapsfrom \max(0,v_{\mathrm{stay}}-1)$ and $v_{\mathrm{change}} \mapsfrom \max(0,v_{\mathrm{change}}-1)$.
			\STATE \Indm \textbf{End For}
			\STATE \textbf{For} each lead/lag vehicle pair (A and B) from downstream most \textbf{do}:
			\STATE \Indp Calculate $\mathbf{A}$, $\mathbf{B}$. 
			\STATE \textbf{If} A and B are TVs \textbf{then}
			\STATE \Indp \underline{Play TU game}: Find $(i^*,j^*)$ and $\sigma$.
			\STATE \Indm \textbf{Else}
			\STATE \Indp \underline{Play NTU game}: Find $(i^{\prime},j^{\prime})$.
			\STATE \Indm \textbf{End If}
			\STATE Update $v_{\mathrm{stay}}^{\mathrm{A}}$, $v_{\mathrm{change}}^{\mathrm{A}}$, $v_{\mathrm{stay}}^{\mathrm{B}}$, and $v_{\mathrm{change}}^{\mathrm{B}}$.
			\STATE \Indm \textbf{End For}
			\STATE \textbf{For} each vehicle \textbf{do}:
			\STATE \Indp \underline{Update system state}: 
			\STATE \textbf{If} $v_{\mathrm{change}} > v_{\mathrm{stay}}$ \textbf{then}
			\STATE \Indp change lane state
			\STATE \Indm \textbf{End If}
			\STATE $v \mapsfrom \max(v_{\mathrm{change}}, v_{\mathrm{stay}})$.
			\STATE \Indm \textbf{End For}
			\STATE \Indm \textbf{End For}
		\end{algorithmic}
	\end{algorithm}
}
\medskip

\noindent \autoref{A:1} considers a two-lane road represented as a two-dimensional uniform lattice $\mathcal{L} \times \mathcal{I}$, where $\mathcal{L}$ is the longitudinal dimension and $\mathcal{I}$ is cross-sectional dimension (i.e., the lanes).  The spatial discretization is such that each site in the lattice can be occupied by at most one vehicle.  The state of each (occupied) site during a discrete time step is specified by discretized speeds, which can take integer values between 0 and $v_{\max}$, where $v_{\max}$ is the maximum number of cells that can be traversed by a vehicle on the road during a single discrete time step. The simulation procedure is summarized in \autoref{A:1}.  

\section{Experiments}
\label{S:simulation}

\subsection{Numerical example}
\label{SS:example}

Consider two vehicle classes on the road, one with high VOT and one with low VOT.  The equilibrium mean speeds are 38 km/h for TVs with high VOT and 31 km/h for TVs with low VOT. We assume that vehicle A has low VOT and wants to make a lane change; it has two options: change lane and achieve a speed of $v^{\mathrm{A}}_1 = 55$ km/h, or stay and maintain a speed of $v^{\mathrm{A}}_2 = 25$ km/h. vehicle B, with high VOT, is the competing lag vehicle and it has two options: do not give way and achieve a speed of $v^{\mathrm{B}}_1 = 52$ km/h and give way to achieve a speed of $v^{\mathrm{B}}_2 = 45$ km/h.  Following \citep{hossan2016investigating}, the VOT coefficient for a low VOT vehicle is set to $c_{\vot}^{\mathrm{A}} = $10 dollars/h and that for the high VOT vehicle class is set to $c_{\vot}^{\mathrm{B}} = $25 dollars/h. We assumed that they are honest in their VOT. We set $t_a = 3 $ s, $a_1^{\mathrm{A}} = -4$ ${\mathrm{m}}/{\mathrm{s}}^2$, $a_2^{\mathrm{A}} = 1$ ${\mathrm{m}}/{\mathrm{s}}^2$, $a_1^{\mathrm{B}} = -3$ ${\mathrm{m}}/{\mathrm{s}}^2$, and $a_2^{\mathrm{B}} = -1$ ${\mathrm{m}}/{\mathrm{s}}^2$. 
Applying \eqref{E:td}, we have that $t_{\dd}^{\mathrm{A}} = 2.26$ s and $t_{\dd}^{\mathrm{B}} = 0.34$ s. Then based on \autoref{t2} and Equation \eqref{E:t3}, we have the payoff matrix given in \autoref{t4}.
\begin{table}[h!]
	\centering
	\caption{Utility matrix for example TU game}
	\begin{tabular}{|c|c|c|c|}
		\hline
		\multicolumn{2}{|c|}{\multirow{2}{*}{Actions}} &
		\multicolumn{2}{|c|}{vehicle B} \\
		\cline{3-4}
		\multicolumn{2}{|c|}{} & Do not give way & Give way \\
		\hline
		\multirow{2}{*}{vehicle A}
		& Change lane & ($-M$, $-M$) & ($0.0062$, $0$) \\
		\cline{2-4}
		& Do not change lane & ($0$, $0.0023$) & ($0$, $0$) \\
		\hline
	\end{tabular}
	\label{t4}
\end{table}

Following the steps in \autoref{SS:TUgame}, we get $\omega^* = 0.0062$, $T_{\mathrm{A}} - T_{\mathrm{B}} = 0$, and the TU solution for this game is $(0.0031, 0.0031)$. The cooperative strategy $(i^*,j^*)$ is \{A changes lanes, B gives way\}, and $A_{i^*j^*} = A_{12} = 0.0062$. Since $\sigma = A_{i^*j^*} - \frac{T_{\mathrm{A}} - T_{\mathrm{B}} + \omega^*}{2} = 0.0062 - 0.0031 = 0.0031$ is positive, A would pay B 0.0031 dollars to change lanes and B receives 0.0031 dollars to give way to A.  It is worth noting that, in this case, even though A has a lower VOT than B, it is still possible that A is willing to pay B to change lanes. But, in general, vehicles with higher VOT are more likely to be the payers.

\subsection{Simulation experiments: Analysis of benefit to TVs}
\label{SS:cost-effectiveness}

\begin{figure}[h!]
	\centering
	\subfloat[][high VOT TVs]{
		\includegraphics[width=0.5\textwidth]{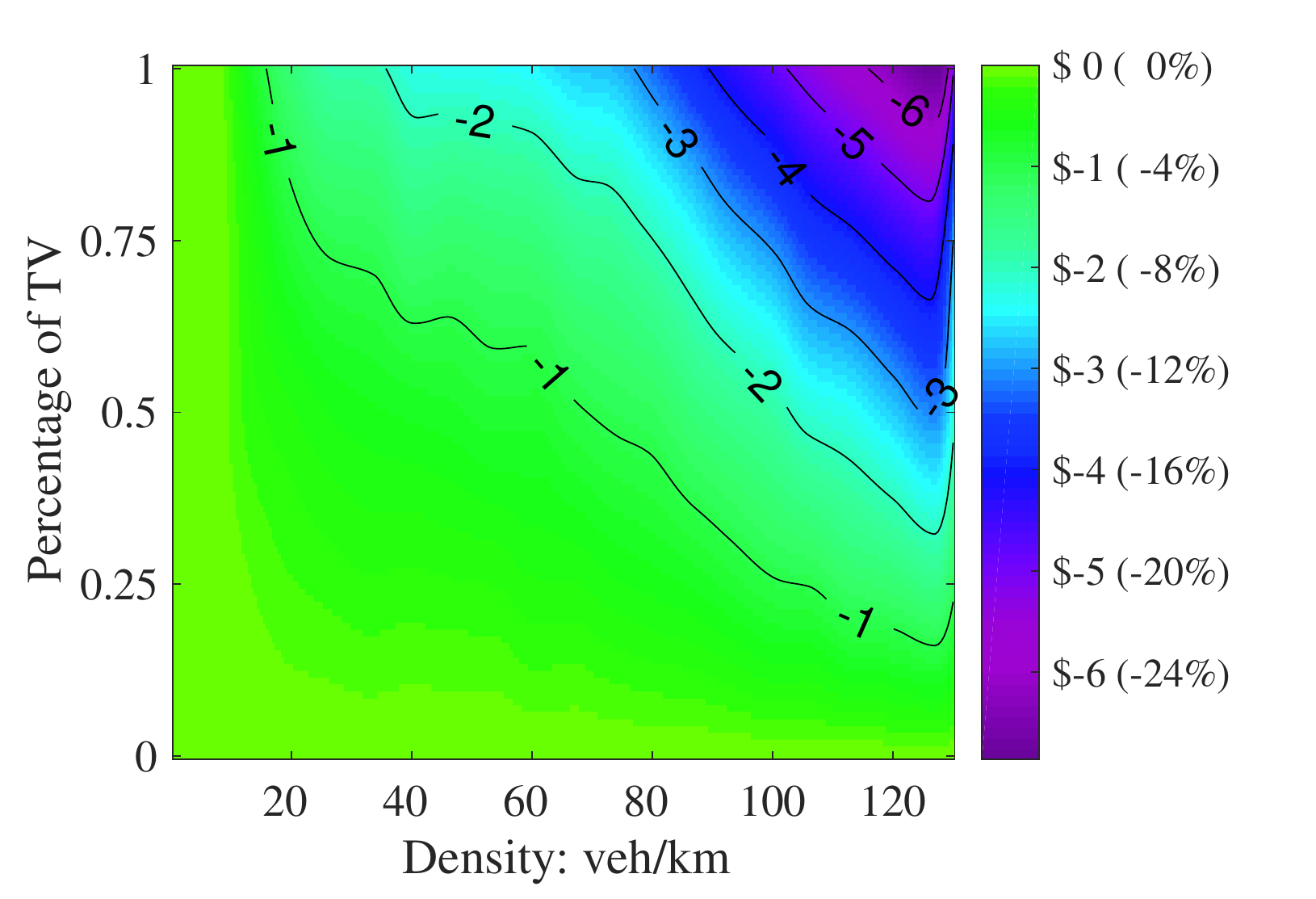}
		\label{F:heat_incomea}} 
	\subfloat[][low VOT TVs]{
		\includegraphics[width=0.5\textwidth]{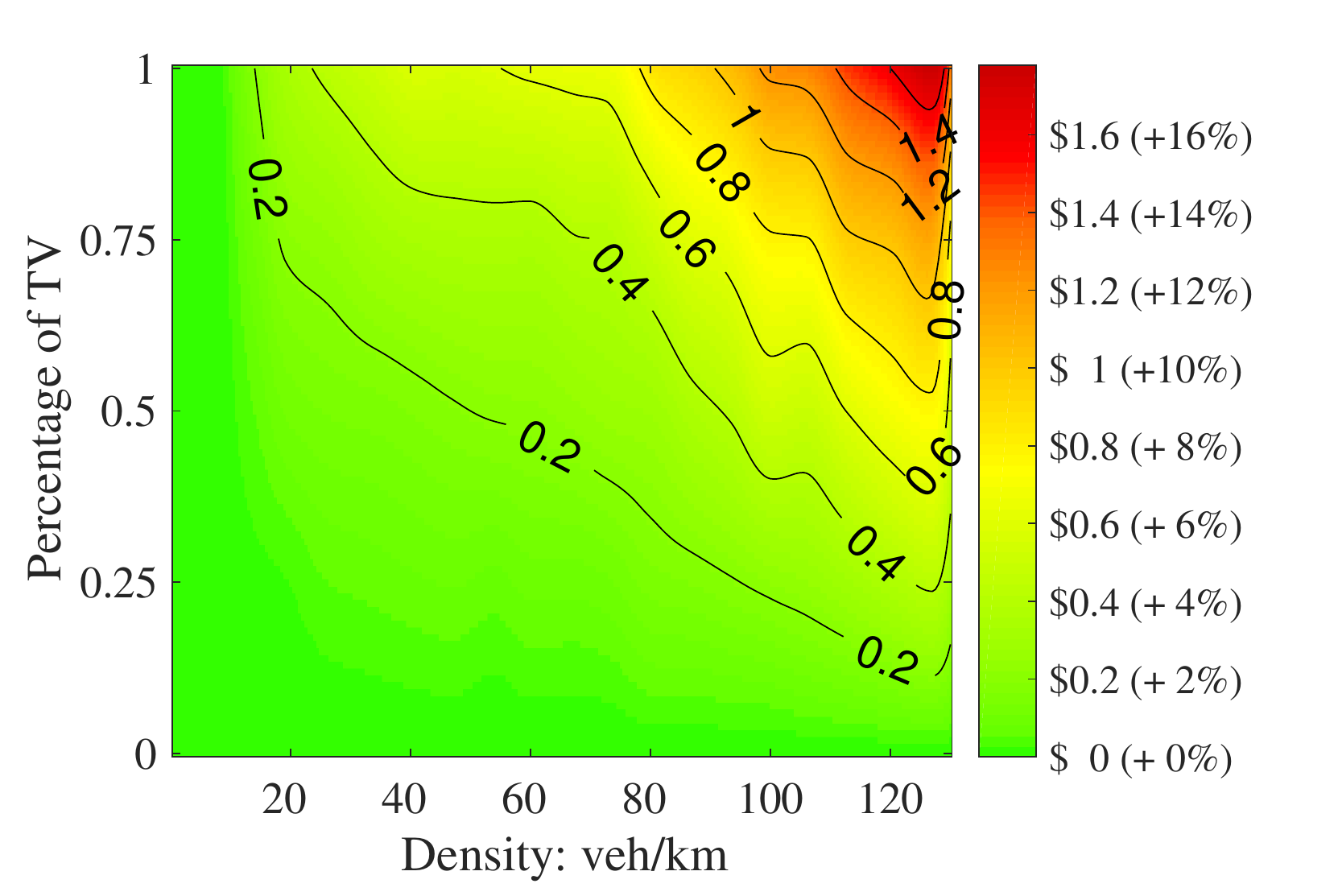}
		\label{F:heat_incomeb}}
	\caption{Income per hour-travel for (a) high VOT TVs and (b) low VOT TVs.}
	\label{F:heat_income}
\end{figure}

\begin{figure}[h!]
	\centering
	\subfloat[][high VOT TVs]{
		\includegraphics[width=0.5\textwidth]{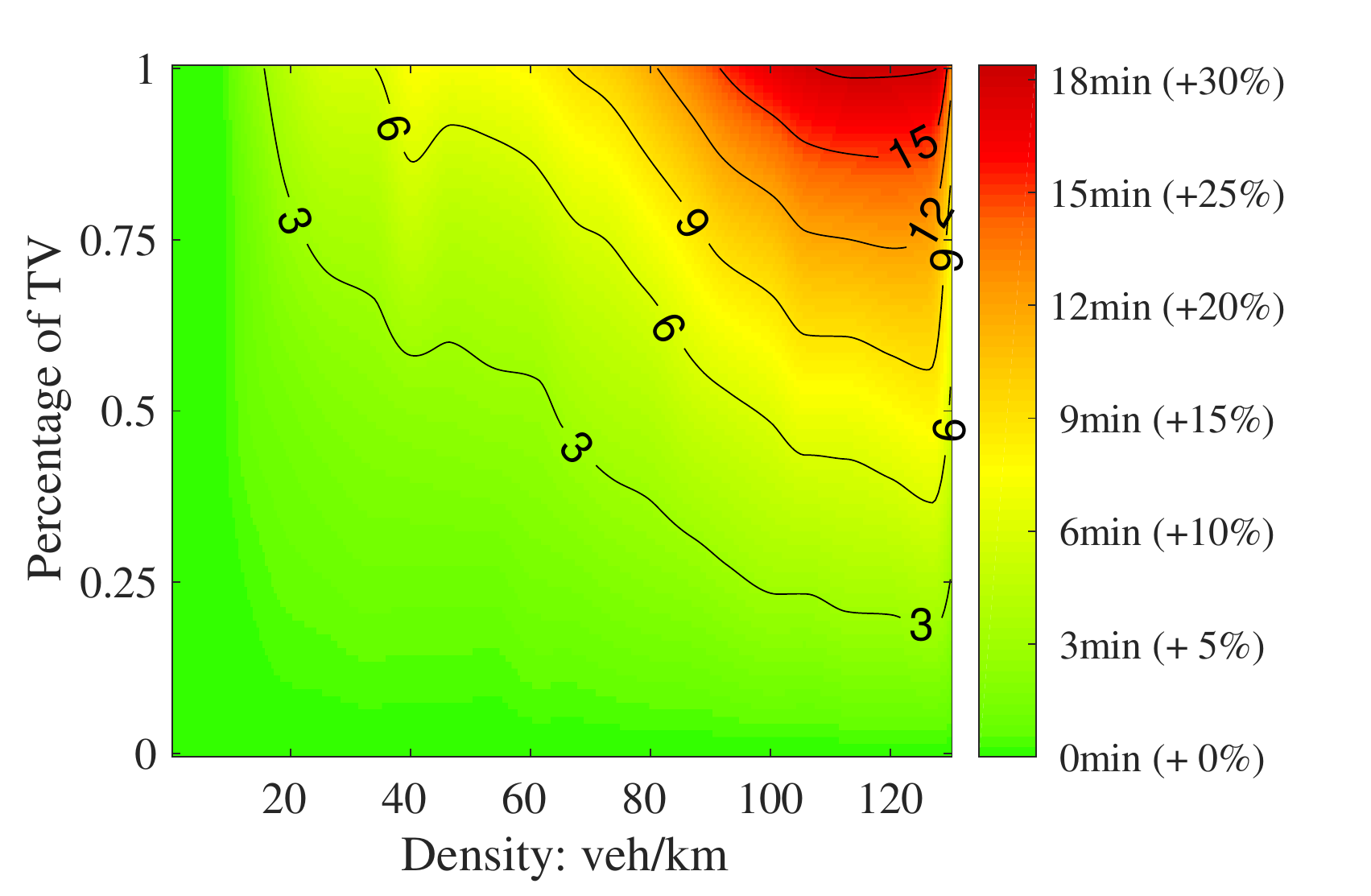}
		\label{F:heat_timea}} 
	\subfloat[][low VOT TVs]{
		\includegraphics[width=0.5\textwidth]{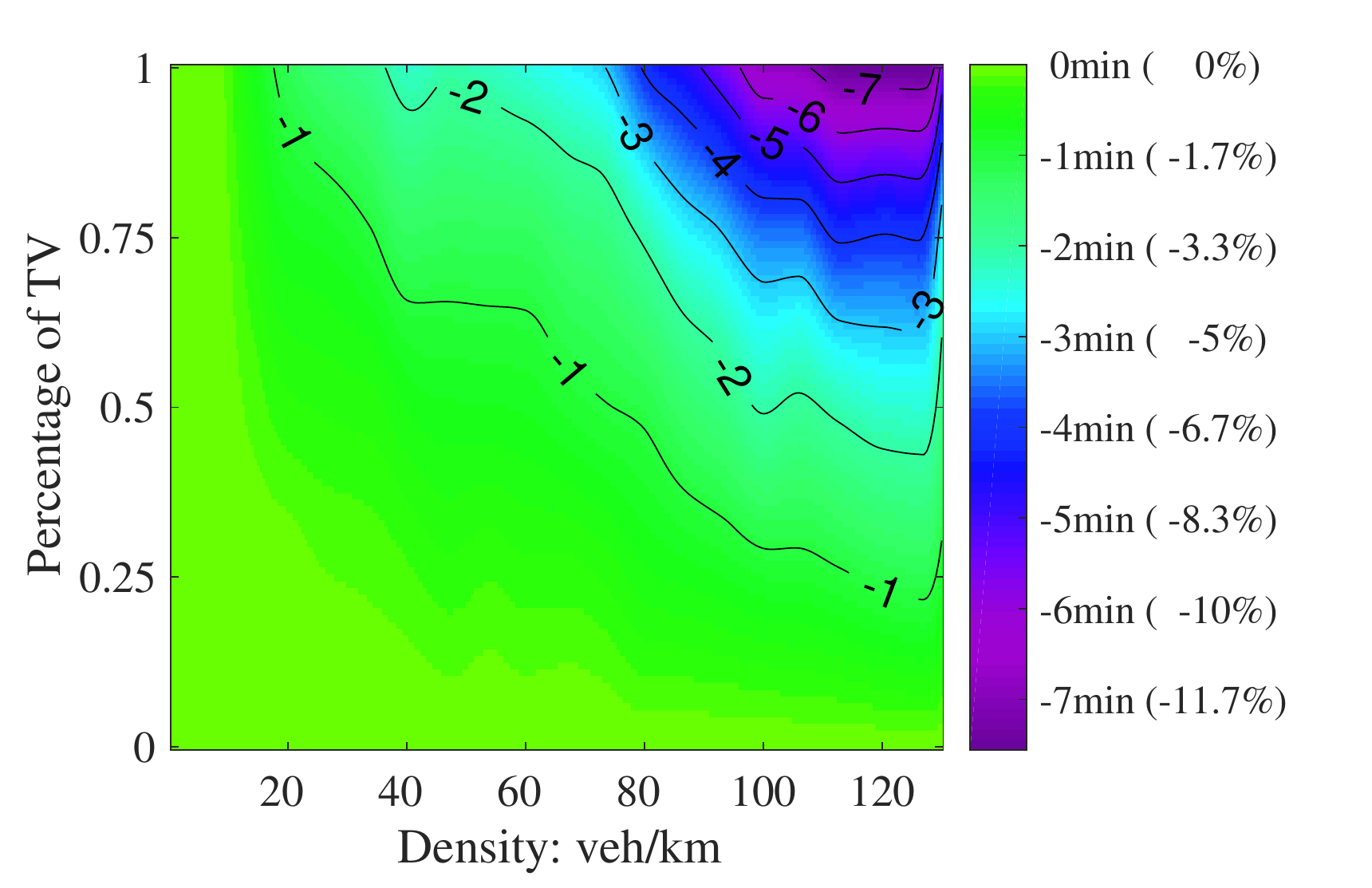}
		\label{F:heat_timeb}}
	\caption{The time saving per hour-travel for (a) high VOT TVs and (b) low VOT TVs.}
	\label{F:heat_time}
\end{figure}

\begin{figure}[h!]
	\centering
	\subfloat[][high VOT TVs]{
		\includegraphics[width=0.5\textwidth]{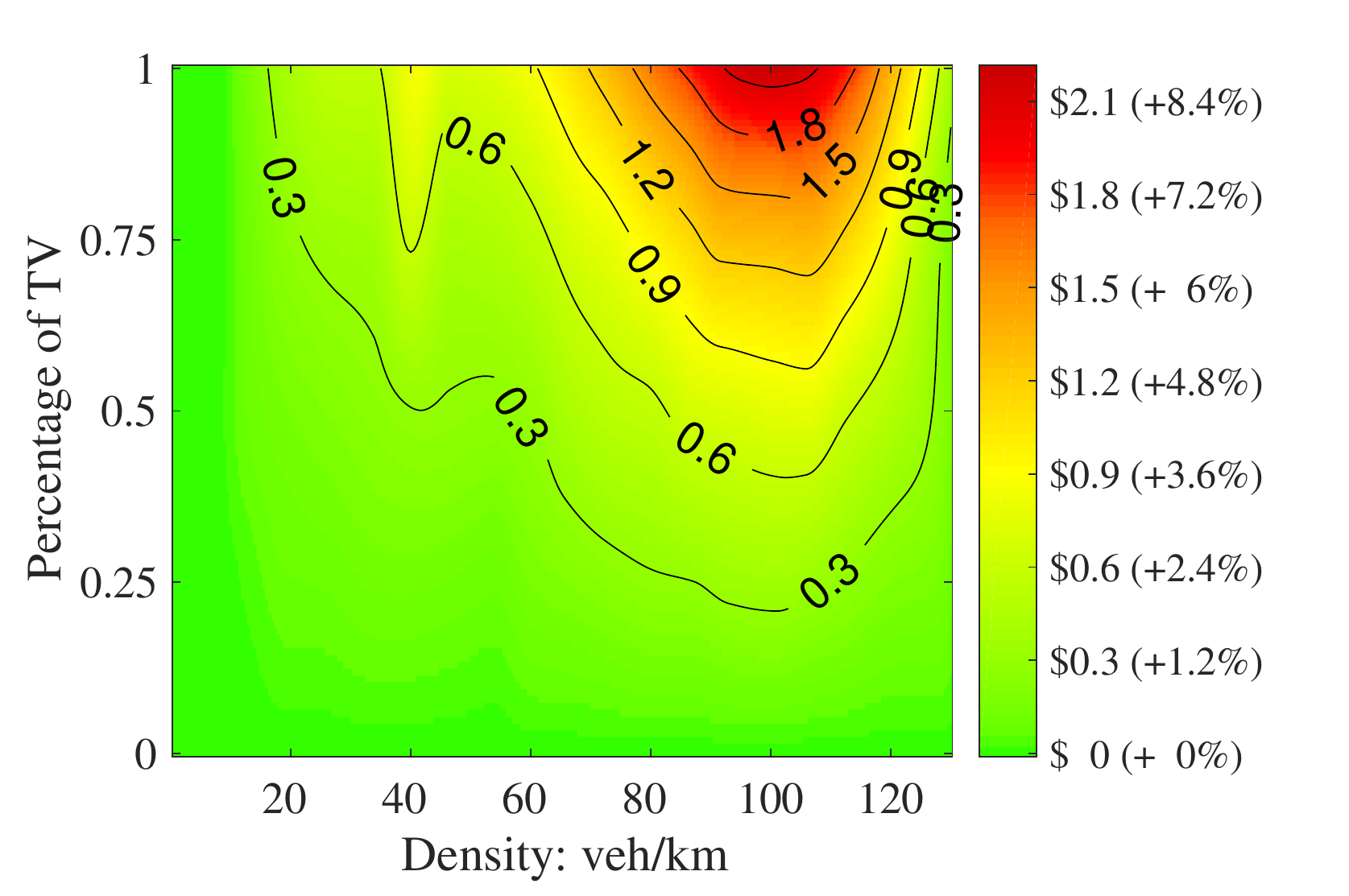}
		\label{F:heat_benefita}} 
	\subfloat[][low VOT TVs]{
		\includegraphics[width=0.5\textwidth]{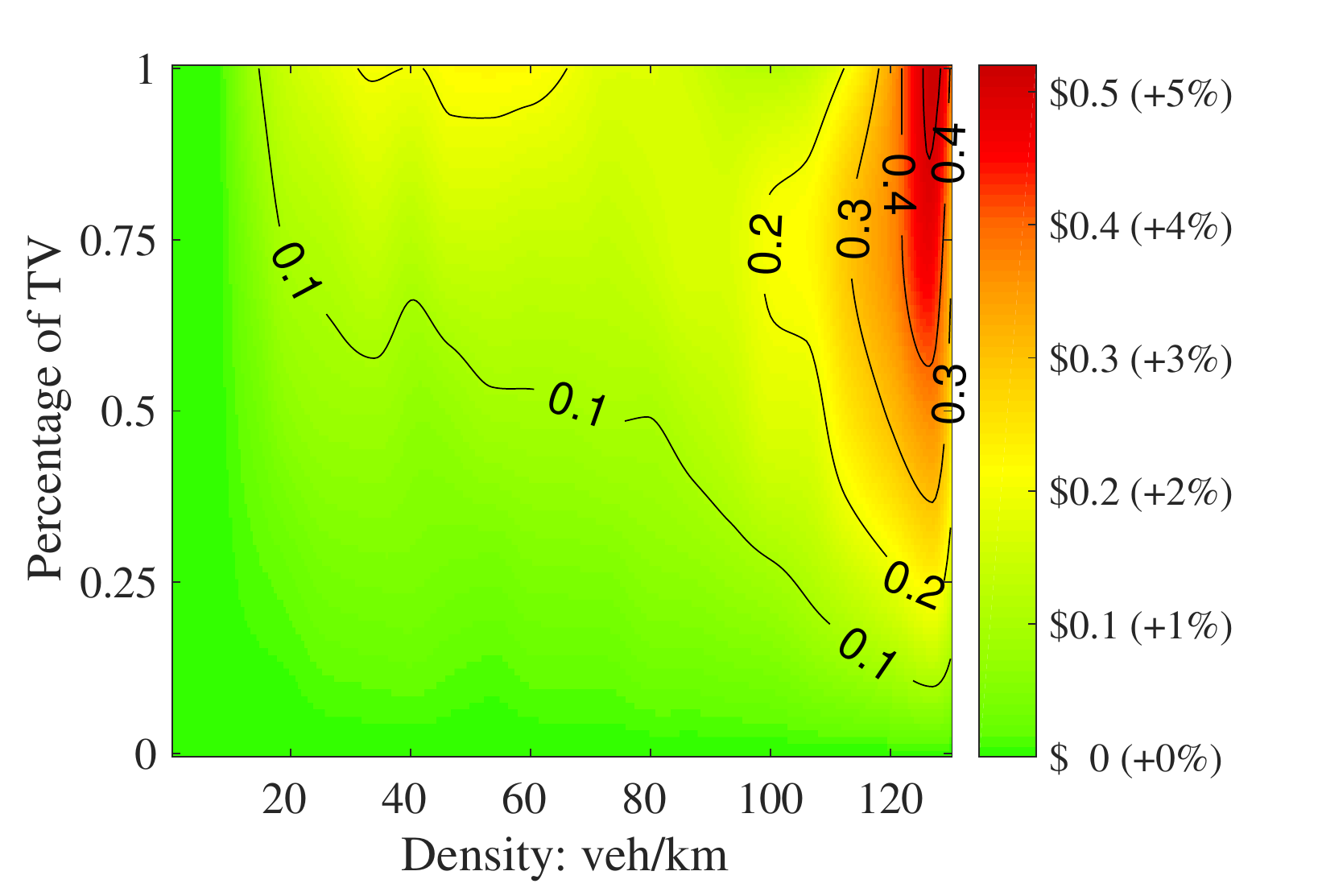}
		\label{F:heat_benefitb}}
	\caption{The total benefit per hour-travel for (a) high VOT TVs and (b) low VOT TVs.}
	\label{F:heat_benefit}
\end{figure}

\begin{figure}[h!]
	\centering
	\subfloat[][untruthful high VOT TVs]{
		\includegraphics[width=0.5\textwidth]{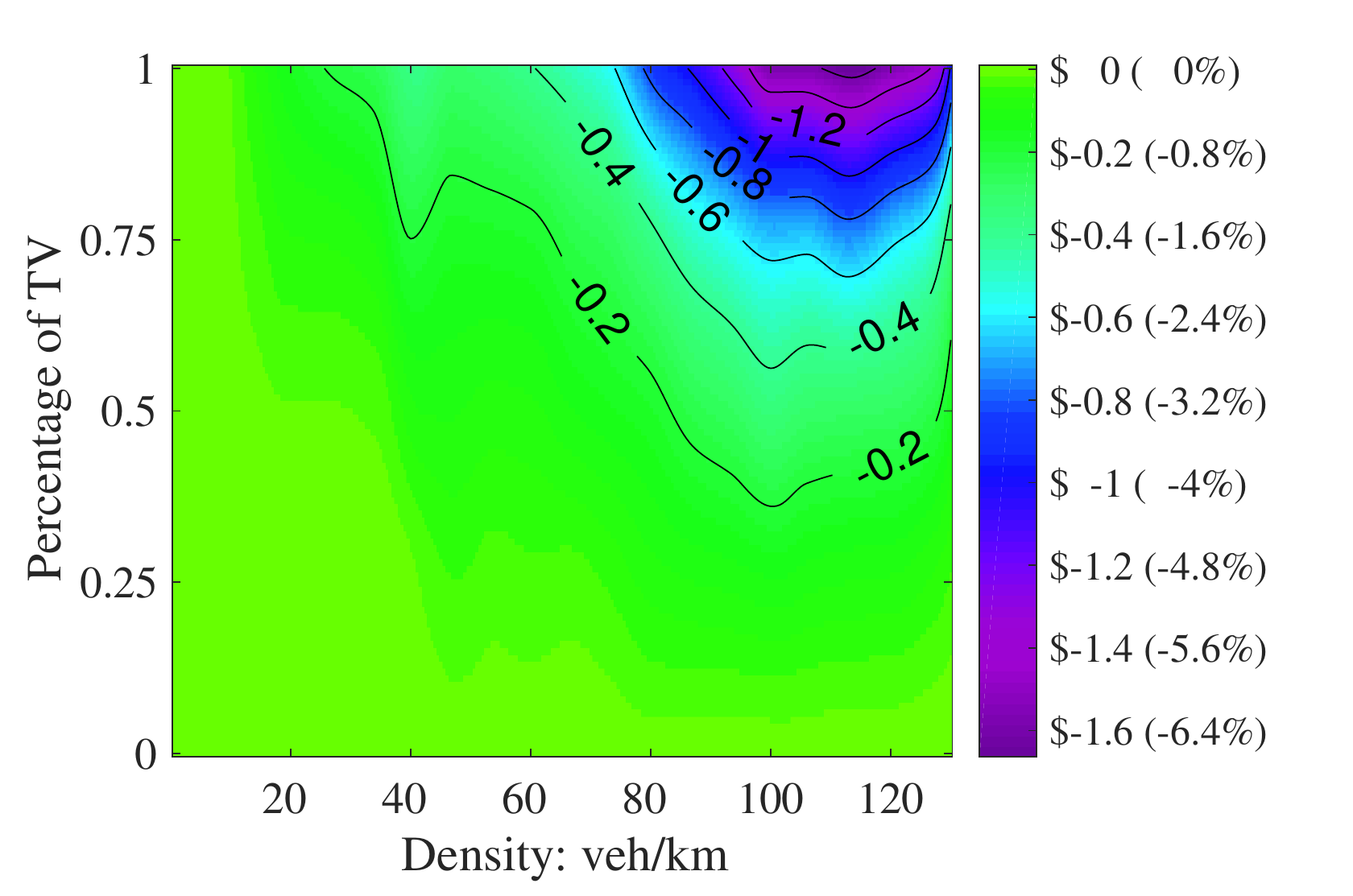}
		\label{F:heat_lyinga}} 
	\subfloat[][untruthful low VOT TVs]{
		\includegraphics[width=0.5\textwidth]{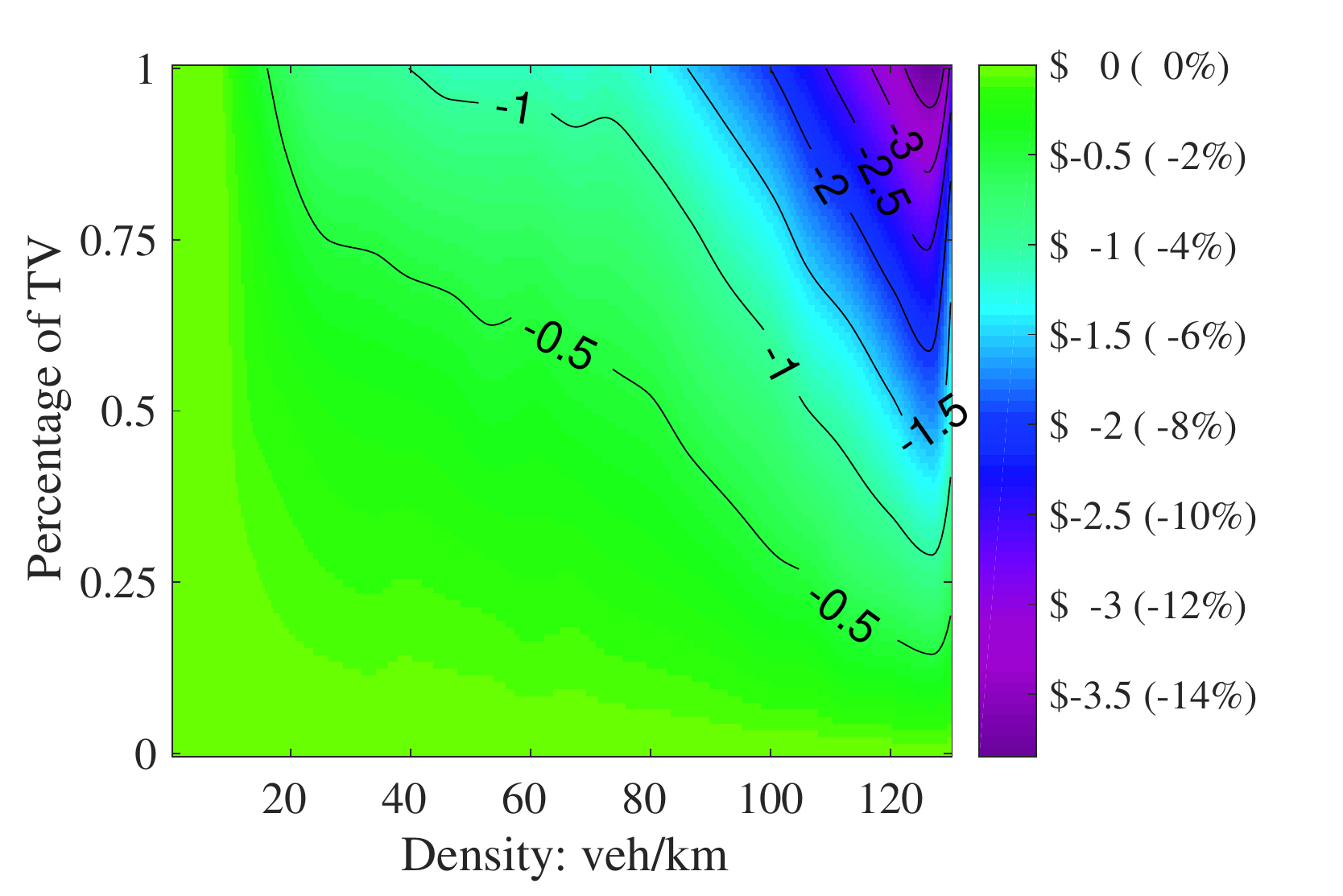}
		\label{F:heat_lyingb}}
	\caption{Benefit per hour-travel for (a) high VOT TVs and (b) low VOT TVs, when they are untruthful about their reported VOT.}
	\label{F:heat_lying}
\end{figure}

The defaults parameters used for the experiments in the remainder of this section are set to: cell size = 7.5 m, total length of road  = 20.25km, time step = 1 s, $v_{\max}$ = 5 cells/s, $a_{\mathrm{pos}}$ = 1 cell/s$^{\mathrm{2}}$, $a_{\mathrm{neg}}$ = -1 cell/s$^{\mathrm{2}}$, number of lanes = 2, $p_{\mathrm{sd}}$ = 1/3, High-VOT = 25 dollars/h, low-VOT = 10 dollars/h \citep{hossan2016investigating}, penetration rate of TV = 1, and high-to-low VOT ratio = 1:4. 
Figures \ref{F:heat_income} - \ref{F:heat_lying} depict the results of simulation experiments, where we change the penetration rates of TVs and traffic densities.  Income and time savings per hour are illustrated in \autoref{F:heat_income} and \autoref{F:heat_time} under varying traffic densities for both high VOT TVs and low VOT TVs.
%
They show that income is, in general, negative for high VOT TVs while the time saved is positive.  For the low VOT TVs, we see the exact opposite.  This is intuitive.  To appreciate the trade-off between travel time saving and income, consider the following ``benefit'' index:
\begin{equation}
	\beta = \frac{1}{\overline{T}^{\mathrm{tr}}}\sum_{i \in \mathcal{V}}\sum_{g \in \mathcal{G}_i}\big(c_{\vot}^i \Delta T_{i,g}^{\mathrm{tr}} + I_{i,g} \big),
\end{equation}
where $\beta$ denotes benefit, $\overline{T}^{\mathrm{tr}}$ is the average travel time of all vehicles $i$ in the set of TV vehicles $\mathcal{V}$, $\Delta T_{i,g}^{\mathrm{tr}}$ is the travel time saved (can be negative) by vehicle $i$ in the lane change game $g$, $\mathcal{G}_i$ is the set of lane change games played by vehicle $i$ and $I_{i,g}$ is the net income earned by vehicle $i$ in game $g$ (can be negative). The total benefit per hour-travel is illustrated in \autoref{F:heat_benefit} under varying densities for high and low VOT TVs.

For both high and low VOT TVs, the benefit is positive in general. Furthermore, as TV penetration rates increase, the total benefit tends to grow.  We see the same pattern as the traffic density increases from 0 to 120 veh/km.  For densities approaching the jam density (133 veh/km in our simulations), lane changing becomes more difficult and benefit approaches zero.  Hence, the highest benefit to both high and low VOT TVs is around heavy congestion, but where vehicles can still perform lane change maneuvers. Similarly, the simulation results (\autoref{F:heat_benefit}) show that in free flow conditions (density $<$10 veh/km), in total jam condition (density $>$131 veh/km) or when TV penetration rate is very low ($<$0.05), the benefit is very small (between -0.2\% and 0.2\%). 
Because the benefit of TVs is always positive, some NTVs may be encouraged to join TU games. Moreover, higher penetration rates mean higher benefit. Hence joining TU games can result in increased benefit to all vehicles.

\textbf{Truthfulness in reporting VOT}. Truthfulness is a topic that has received little attention in traffic flow research \citep{yang2018auction}. To investigate the impact of untruthful reporting of VOT, we relax the truthful reporting assumption and consider two scenarios: one in which high-VOT vehicles declare they are low-VOT vehicles and one where low-VOT vehicles declare they are high-VOT vehicles.  As with the previous experiments, we test these two scenarios under varying traffic densities and varying TV penetration rates. For high VOT vehicles, the experiments attempt to gauge whether the monetary gains (whether this is in the form paying less to change lanes or being paid when denied a lane change opportunity) exceed (on average) the gains in travel time weighted by VOT.  Similarly, in the case of the low VOT vehicle reporting a higher VOT, the experiments are mean to gauge whether the improvement in travel time weighted by true VOT that can be gained by being untruthful can exceed (on average) the amounts they are paid by other vehicle the vehicles.  In both cases, our experiments indicate that the answer in no: being untruthful reduces total benefit, $\beta$. 

\autoref{F:heat_lying} illustrates the total benefits per hour travel for these two scenarios; \autoref{F:heat_lying}\subref{F:heat_lyinga} is a plot of the total benefit of high VOT TVs when they are untruthful (they declare that they are low VOT TVs) and \autoref{F:heat_lying}\subref{F:heat_lyingb} is a plot of the total benefit of low VOT TVs when they are untruthful (they declare that they are high VOT TVs).  In both cases, we observe negative benefits in high traffic density scenarios, which suggests that there is no incentive to lie about their VOT.  This hints that the mechanism that we proposed incentivizes truthfulness.  We leave this at the conjecture level and attack the problem of mechanism design and truthfulness to future research.
%

\begin{figure}[h!]
	\centering
	\resizebox{0.5\textwidth}{!}{%
		\includegraphics{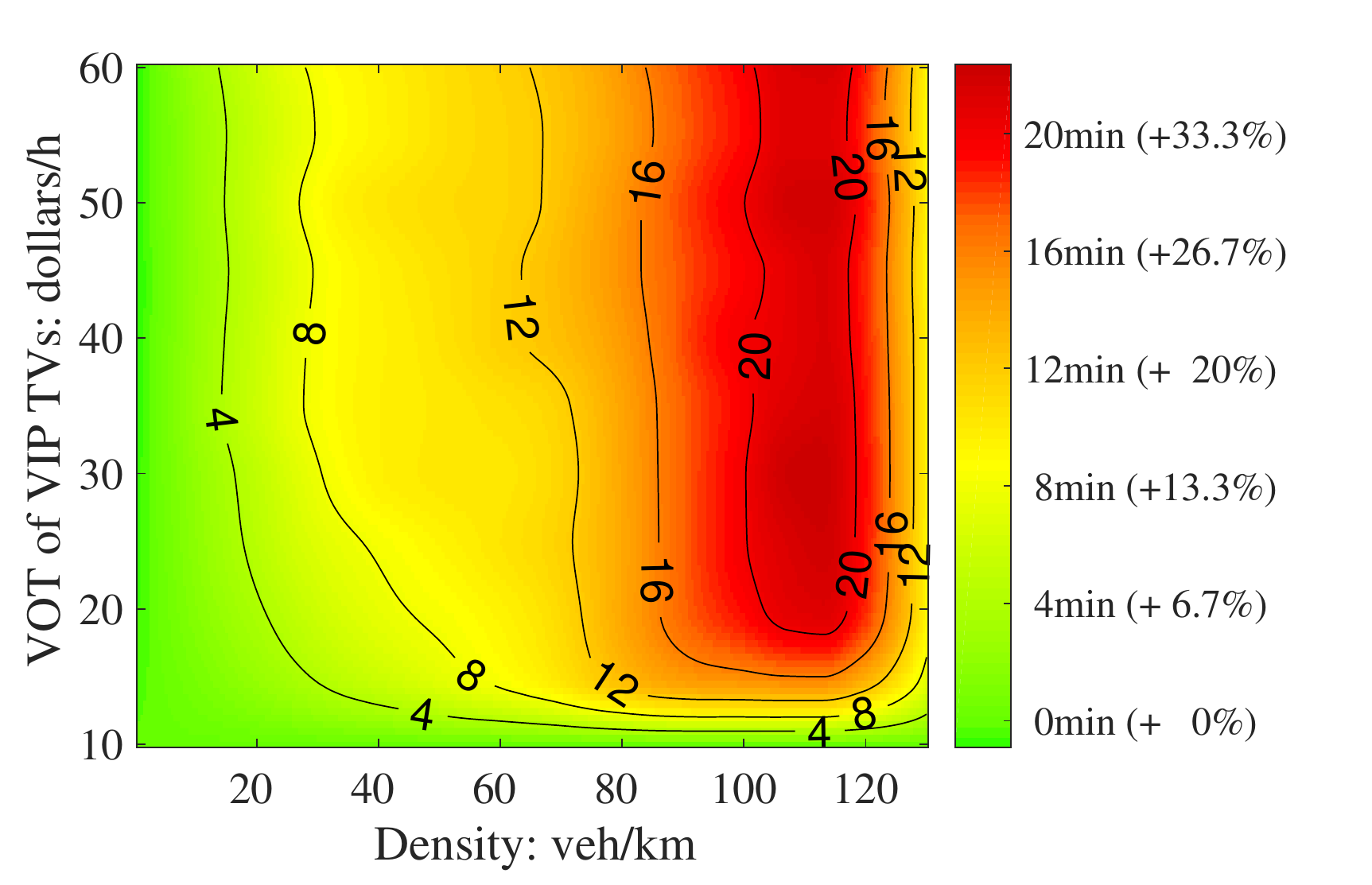}}
	\caption{Saved time per hour-travel when setting different VOTs for high-VOT TVs.} 
	\label{F:saved_hour_high_VOT}
\end{figure}

\textbf{``VIP'' TVs}. In this experiment, we vary the $c_{\vot}$ for the high VOT TV while holding all else fixed. Specifically, $c_{\vot}$ for the low VOT TV is 10 dollars/h and the penetration rate of high VOT TVs is held at 1\% (very small percentage). This is an example of high profile vehicle (along with their entourage).  The results are illustrated in \autoref{F:saved_hour_high_VOT}.
In moderate to very high congestion, the ``VIP'' TVs save about 10-35\% in travel time when $c_{\vot}$ is increased from 10 dollars/h to 60 dollars/h.  Interestingly, however, we observe a bound on time saving of about 38\%, which when reached cannot be improved with greater payment (greater $c_{\vot}$ doesn't make sence when it is larger than 40 dollars/h). This can be attributed to the simple nature of the TU game that is being tested in this paper.  To be specific, the 38\% bound can be attributed to the ``VIP'' TV being blocked by their leaders (in the subject lane).  The 38\% bound may be broken if the TV is allowed to engage in transactions with multiple vehicles simultaneously, namely, including the lag vehicle in the target lane and leaders in both the subject and target lanes.

\subsection{Simulation experiments: Impact on traffic}
\label{SS:characteristic}

\begin{figure}[h!]
	\centering
	\subfloat[][All vehicles vs. NTVs]{
		\includegraphics[width=0.45\textwidth]{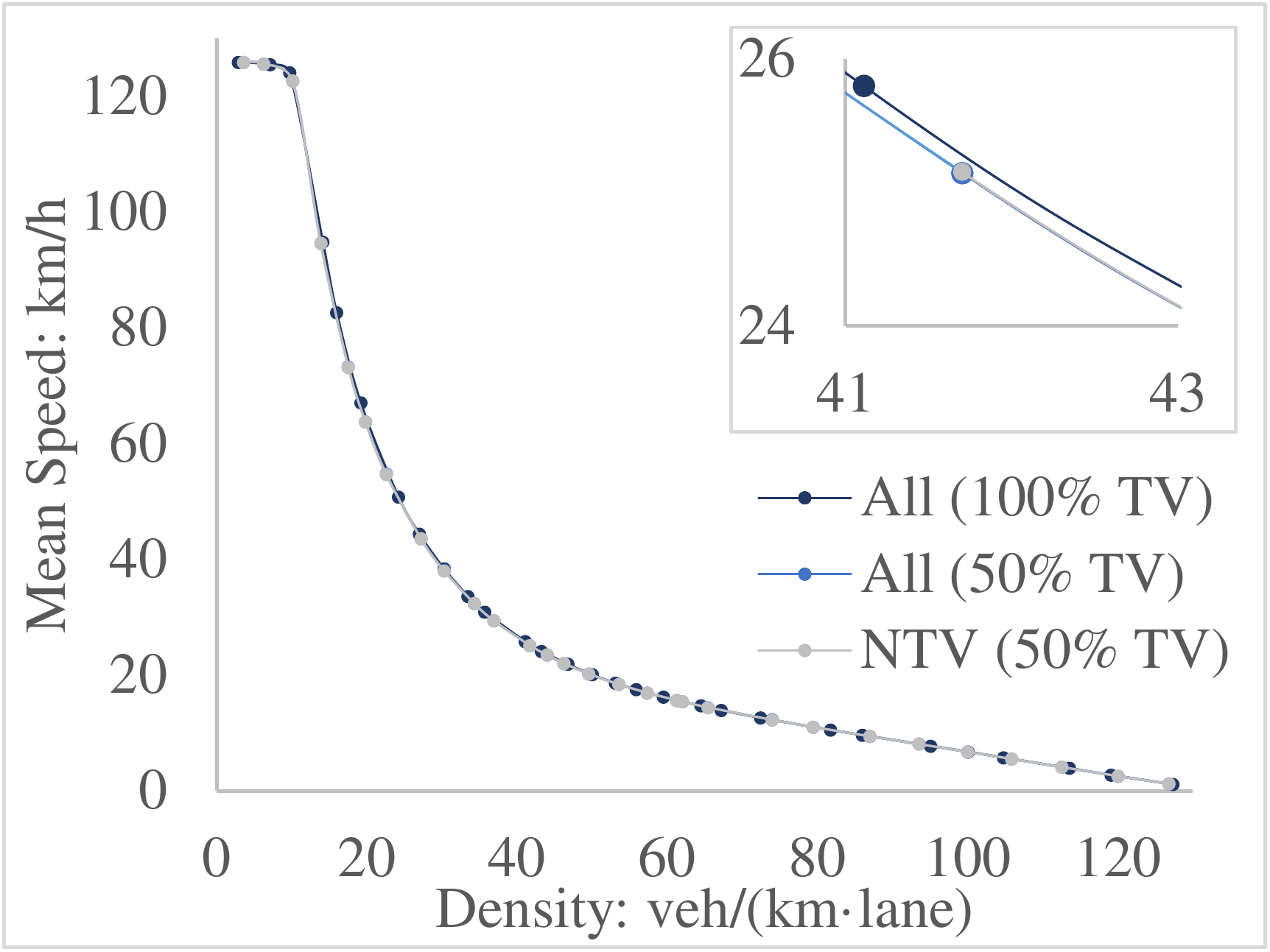}
		\label{F:speed_densitya}} 
	\subfloat[][High-VOT TVs vs. Low-VOT TVs]{
		\includegraphics[width=0.45\textwidth]{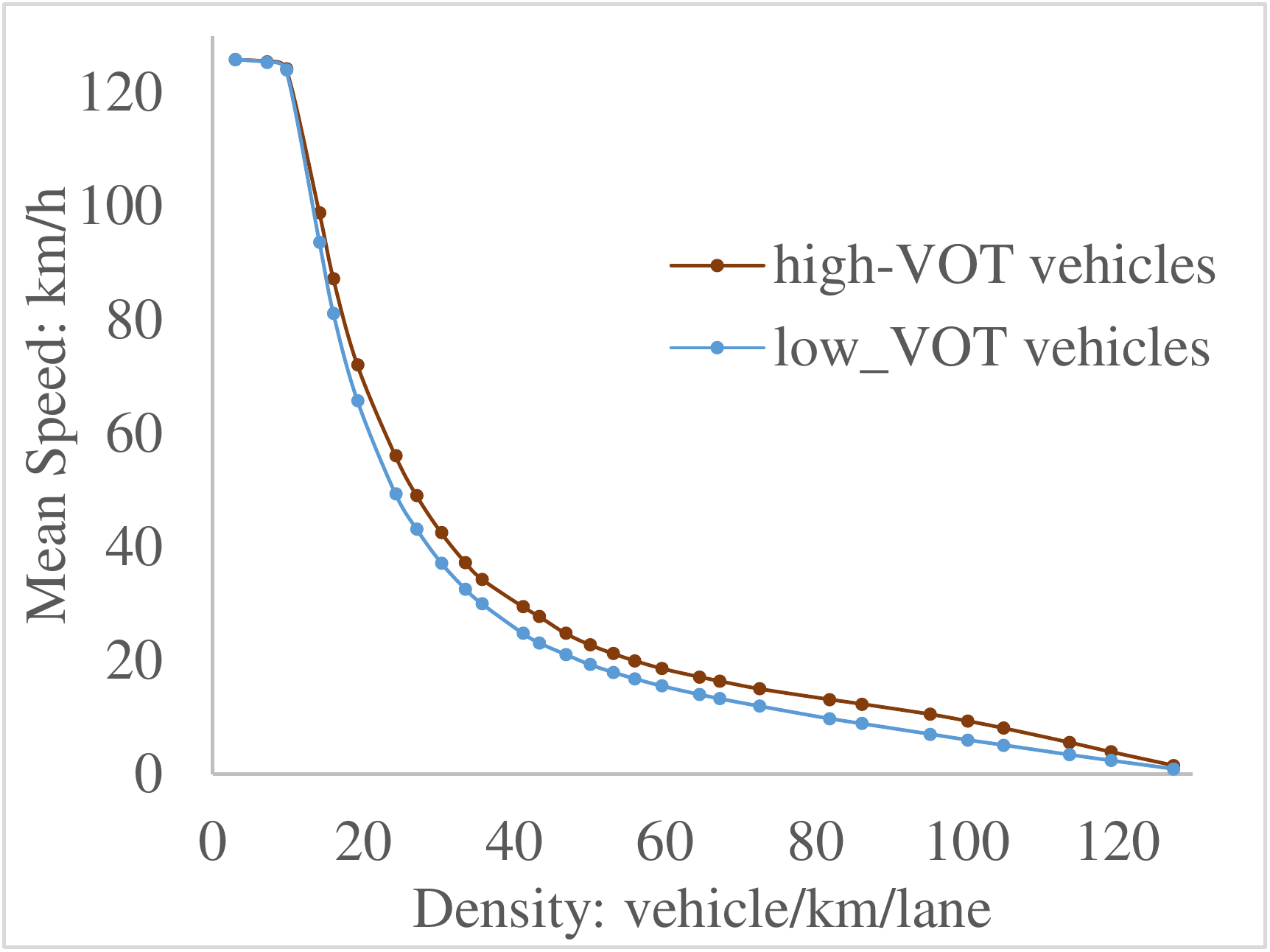}
		\label{F:speed_densityb}}
	\caption{Relationship between global density and mean speed of (a) All vehicles vs. NTVs and (b) High-VOT TVs vs. Low-VOT TVs.}
	\label{F:speed_density}
\end{figure}

\textbf{Speed-density relation}.  We examine the impact of introducing transactions on traffic as whole. The first experiments compare the resulting speed-density relations when TVs constitute 100\% of the vehicle population and when they constitute only 50\% of the vehicle population. \autoref{F:speed_density}\subref{F:speed_densitya} compares the speed-density relations of NTVs in the 50\% case with all vehicles (both TV and NTV) in both the 100\% and 50\% cases.  The figure shows a very small (almost unnoticeable) decrease in mean (equilibrium) speeds of NTVs across a range of traffic densities when the percentage of TVs is 50\%.   We can, at the very least, say that the introduction of TVs does not adversely impact mean speeds (in fact, we see slight improvement). 
In \autoref{F:speed_density}\subref{F:speed_densityb}, we see a more noticeable difference in mean speeds across a range of traffic densities between high and low VOT TVs. Specifically, we see that the speeds of high-VOT TVs are higher than the speeds of low-VOT TVs, outside of free flow conditions and totally jammed traffic, by as much as 60\%.  In the former case, lane changes are either not needed to improve speeds or can be carried out without the need to negotiate gaps; in the latter case, changing lanes will not help improve speeds. 

\begin{figure}[h!]
	\centering
	\subfloat[][high VOT TVs]{
		\includegraphics[width=0.45\textwidth]{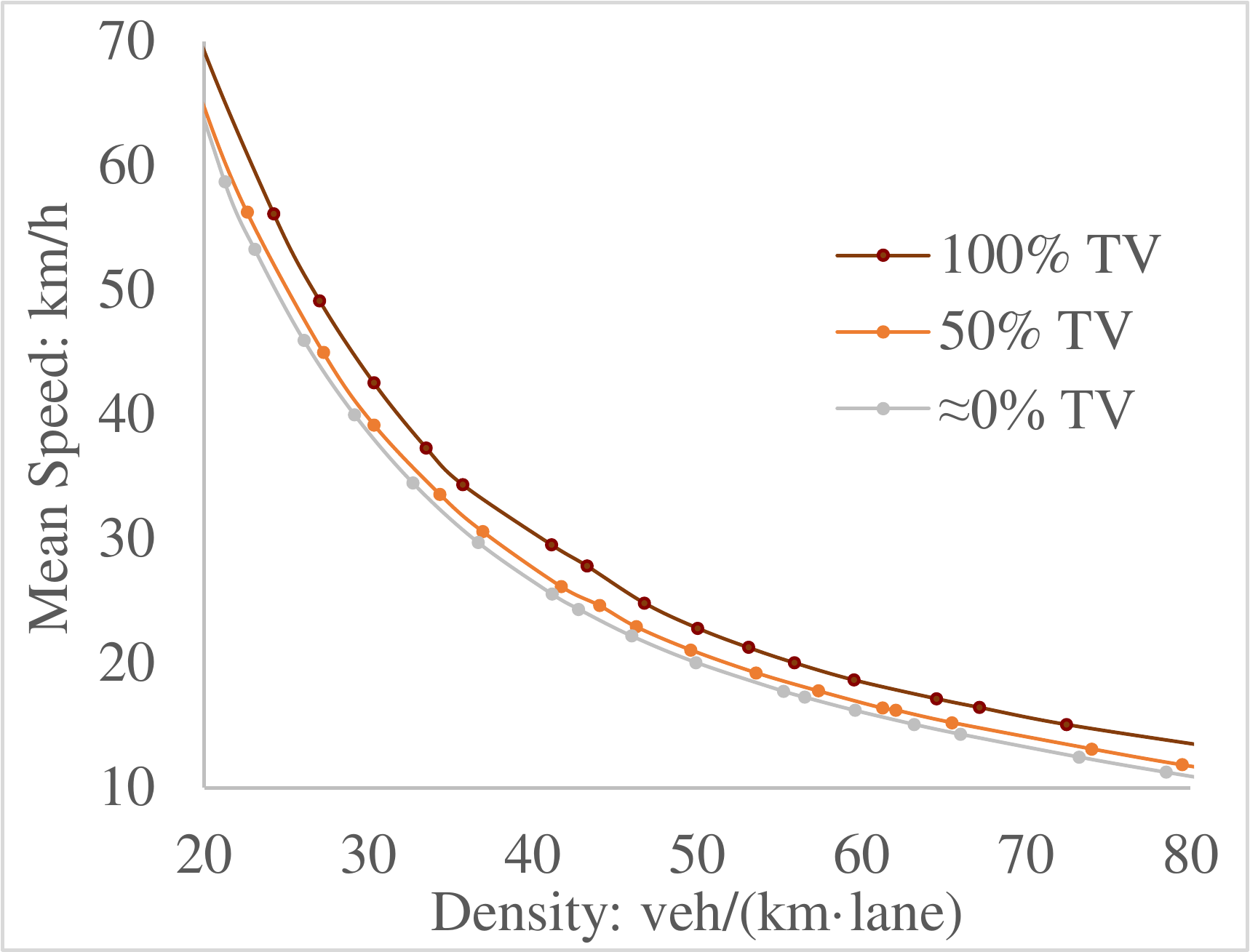}
		\label{F:speed_density_different_penetrationa}} 
	\subfloat[][low VOT TVs]{
		\includegraphics[width=0.45\textwidth]{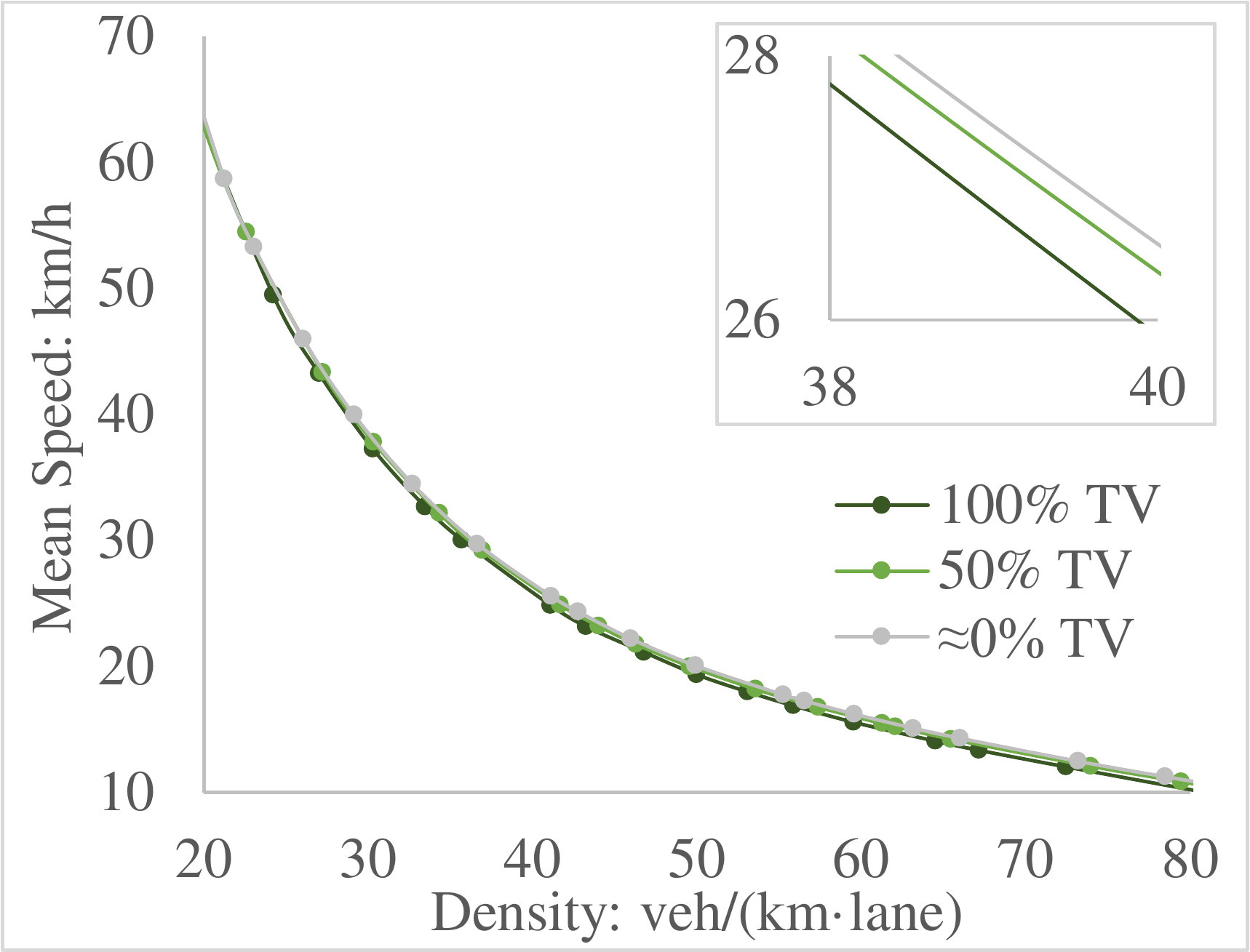}
		\label{F:speed_density_different_penetrationb}}
	\caption{Density and speed relationship varying TVs' penetration rate: (a) high VOT TVs (b) low VOT TVs.}
	\label{F:speed_density_different_penetration}
\end{figure}

We next investigate the impact of TV penetration on the speed-density relations pertaining to high and low VOT TVs.   \autoref{F:speed_density_different_penetration} shows that as the TV penetration rate increases, the speeds of high VOT TVs will increase, the speeds of low VOT TVs tend to drop, but very slightly. (Note that the scale of the $x$-axis in \autoref{F:speed_density_different_penetration}\subref{F:speed_density_different_penetrationa} and \autoref{F:speed_density_different_penetration}\subref{F:speed_density_different_penetrationb} do not include the entire range of traffic densities that were tested; this was done to illustrate the differences/similarities for the different penetration rates.)
As the TV penetration rates approach zero, the speed-density relations of the high and the low VOT TVs converge. The higher TV penetration rate means a higher probability of TU games taking place, and this gives the high VOT TVs more opportunities to increase their speeds, while low VOT TVs can also ``sell" their time more frequently. 
 
\begin{figure}[h!]
	\centering
	\subfloat[][high VOT TVs]{
		\includegraphics[width=0.45\textwidth]{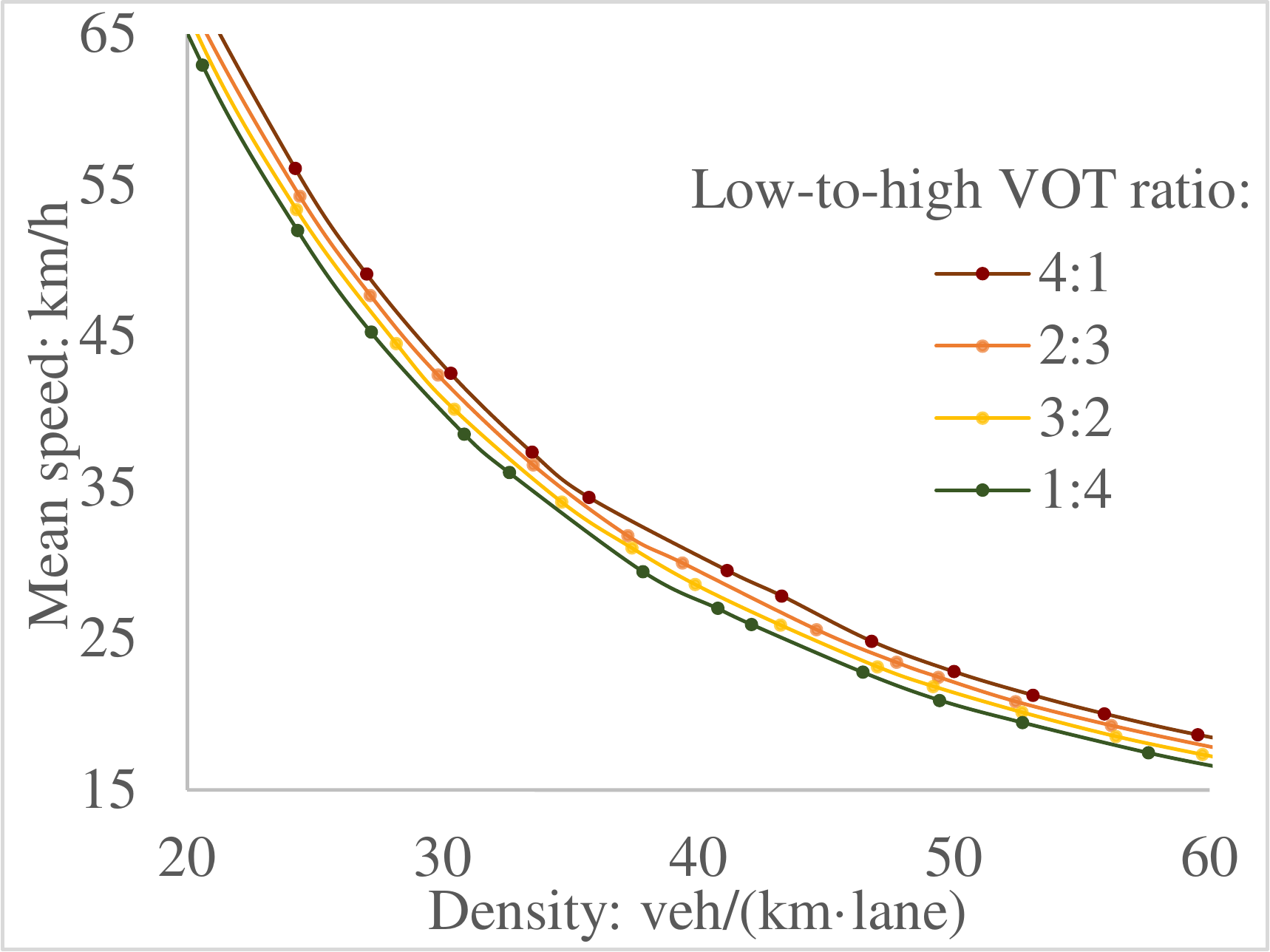}
		\label{F:speed_density_different_ratioa}} 
	\subfloat[][low VOT TVs]{
		\includegraphics[width=0.45\textwidth]{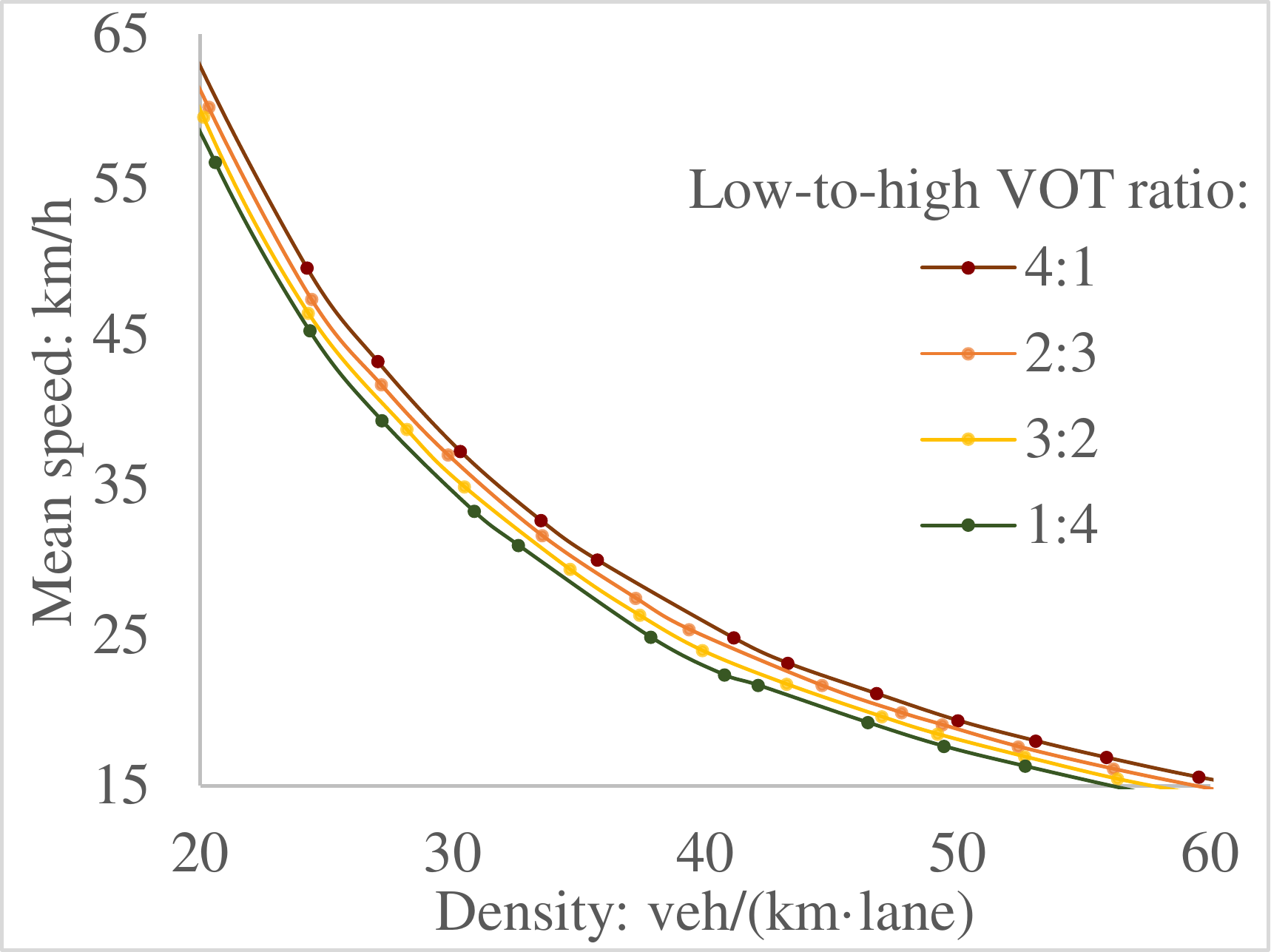}
		\label{F:speed_density_different_ratiob}}
	\caption{Density and speed relationship varying ratio of high VOT TVs and low VOT TVs: (a) high VOT TVs, (b) low VOT TVs.}
	\label{F:speed_density_different_ratio}
\end{figure}

We next investigate the impact of varying low-to-high VOT TV ratios.  The results are given in \autoref{F:speed_density_different_ratio}. Clearly, high VOT TVs always have higher speeds than low VOT TVs. An interesting finding is that, as the ratio of high VOT TVs increases, the speeds of both low and high VOT vehicles decrease. On one hand, a higher high VOT TV fraction results in lower frequencies of transactions with low VOT TVs, leading to a decrease in high VOT TV speeds. On the other hand, higher high VOT TV fractions results in low VOT TV giving way to high VOT TVs with frequency, leading to a decrease in low VOT TV speeds. Hence a healthy market share should have a relatively higher percentage of low VOT vehicles than high VOT vehicles.

\begin{figure}[h!]
	\centering
	\subfloat[][100\% TVs, 13.3 veh/km$\cdot$lane]{
		\includegraphics[width=0.5\textwidth]{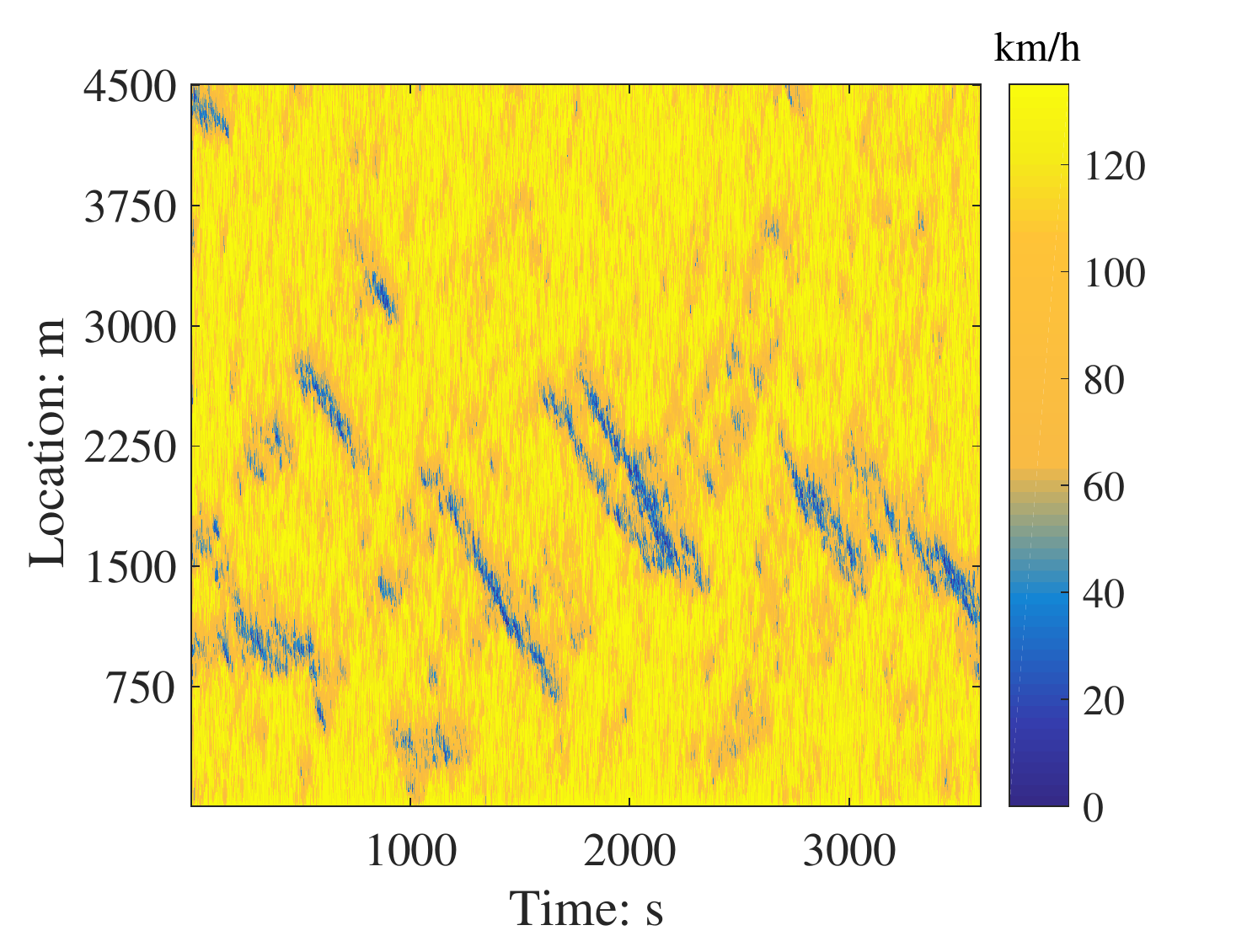}
		\label{F:speed_heat_mapa}} 
	\subfloat[][100\% TVs, 33.3 veh/km$\cdot$lane]{
		\includegraphics[width=0.5\textwidth]{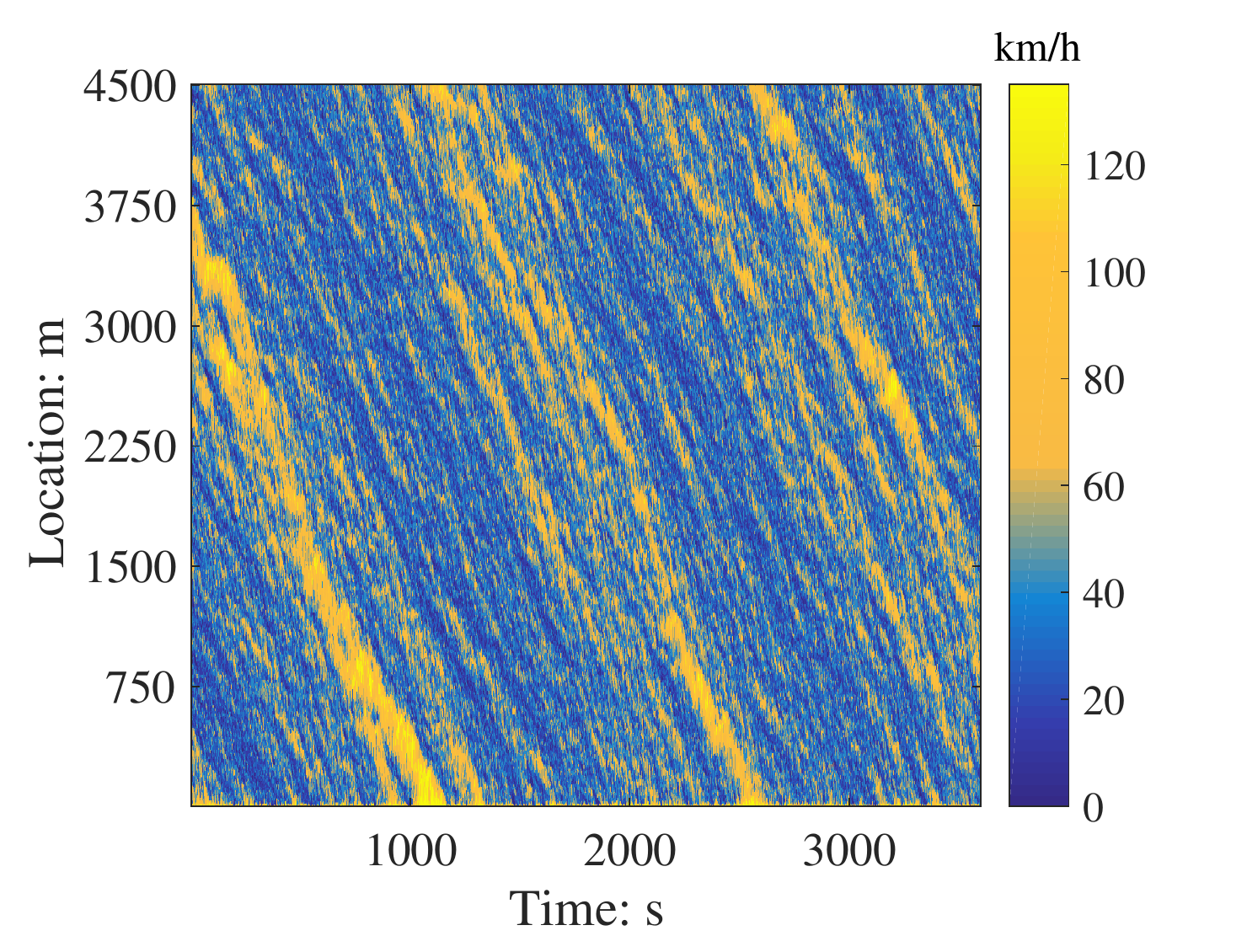}
		\label{F:speed_heat_mapb}}
	
	\subfloat[][100\% NTVs, 13.3 veh/km$\cdot$lane]{
		\includegraphics[width=0.5\textwidth]{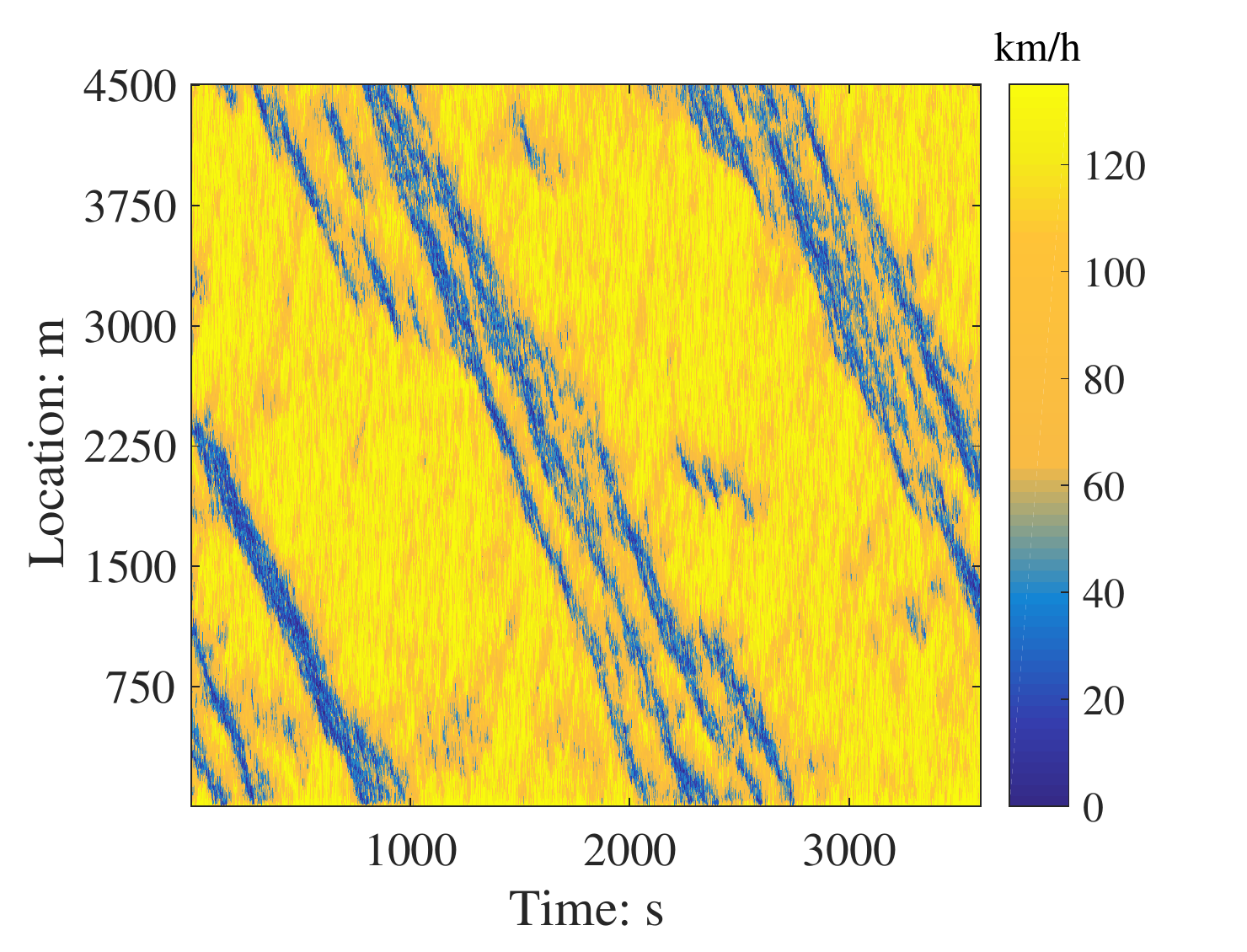}
		\label{F:speed_heat_mapc}}
	\subfloat[][100\% NTVs, 33.3 veh/km$\cdot$lane]{
		\includegraphics[width=0.5\textwidth]{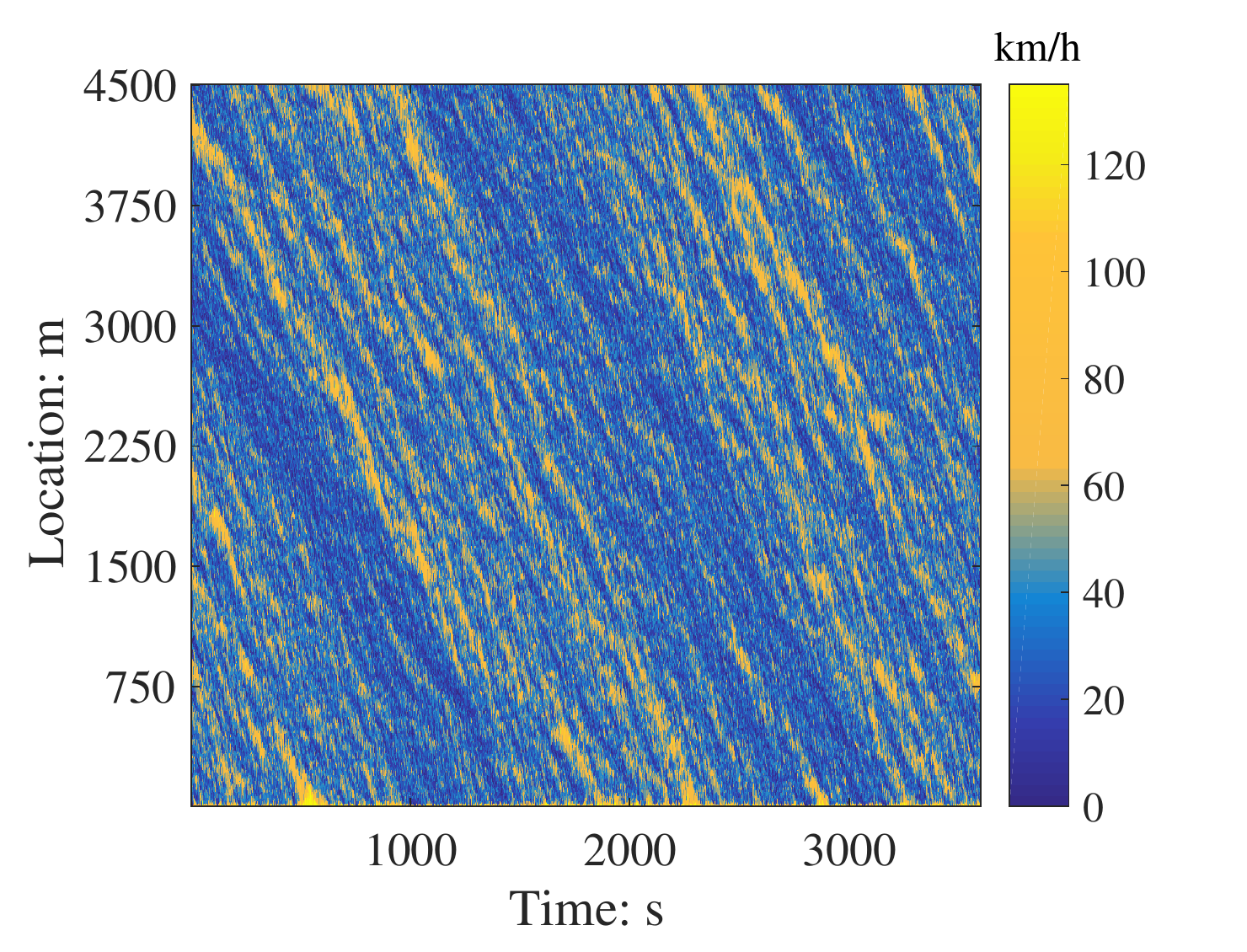}
		\label{F:speed_heat_mapd}}
	\caption{Speed heat maps: (a) 100\% TVs with 13.3 veh/km$\cdot$lane, (b) 100\% TVs with 33.3 veh/km$\cdot$lane, (c) 100\% NTVs with 13.3 veh/km$\cdot$lane, and (d) 100\% NTVs with 33.3 veh/km$\cdot$lane.}
	\label{F:speed_heat_map}
\end{figure}

\textbf{Shock wave formation}. Finally, a 4.5km ring road is simulated over a 1-hour time period and we compare the traffic dynamics that arise in two scenarios: (i) 100\% TVs (i.e., 0\% NTVs) and (ii) 0\% TV (i.e., 100\% NTVs).  Two cases are investigated in each scenario: a case of below (but near) critical traffic density with 13.3 veh/km$\cdot$lane and a case of super-critical (jammed) traffic with 33.3 veh/km$\cdot$lane.    \autoref{F:speed_heat_map} shows the resulting speed heat maps (average of both lanes) for all four cases. \autoref{F:speed_heat_map}\subref{F:speed_heat_mapa} and \autoref{F:speed_heat_map}\subref{F:speed_heat_mapc} are the heatmaps obtained in the sub-critical cases (average density = 13.3 veh/km$\cdot$lane).  We see the formation of both forward and backward waves clearly in both cases, but they appear to be less severe in the first scenario (100\% TVs). \autoref{F:speed_heat_map}\subref{F:speed_heat_mapb} and \autoref{F:speed_heat_map}\subref{F:speed_heat_mapd} are the heatmaps obtained in the super-critical cases (average density = 33.3 veh/km$\cdot$lane).  
In general, there is no significant difference in the wave formation characteristics in high density traffic conditions. Therefore, while it is difficult to conclude that TU games can prevent stop and go waves from forming, we can comfortably conclude that TU games do not have an adverse effect on traffic conditions in general.

\section{Conclusion and outlook}
\label{S:Conc}
Vehicles available on the market today come equipped with advanced sensor technologies and in many vehicles features such as adaptive cruise control and lane departure warning come standard. This allows vehicles (on the roads today) to respond to traffic conditions around them.  It is safe to bet that communication between vehicles is right around the corner.  This creates opportunities for re-thinking traffic management in radical ways. This thinking is the premise of this paper.  We proposed a treatment of lane changing as transferable utility (TU) games with side payments.  We demonstrated that the proposed utility transfer allows vehicles engaged in such transactions to achieve pareto efficient payoffs.  The main idea is that vehicles exchange right-of-way for money.  This constitutes a departure from cooperative lane changing involving winners \textit{and} losers to one where all players can win.

A cellular automaton was developed to perform experiments.  The simulation results indicate that the ability to play TU games had no impact on travel times in free-flow conditions, heavy (bumper-to-bumper) traffic conditions, or when the penetration of vehicles willing to engage in TU games approaches zero.  Otherwise, both vehicles with low and high values of travel time (VOT) derived benefit from the proposed approach.

The proposed model is rather simple. Our aim was to test the idea of TU games for lane changing at an approachable level.  This creates numerous avenues for future research.  From considering more sophisticated models of utility and traffic dynamics, to utility transfer for mandatory lane change maneuvers, to games that involve more than two players with more choices, to games that include unconnected vehicles.  Along the lines of the utility function, the estimated values of the constants can have a significant impact on the time difference parameter, which can have a substantial impact on the outcomes of the games.  One topic for future research involves the design of lane changing games that are robust to estimation errors.  On the other hand, in the absence of estimation errors, we conjecture that the TU framework proposed in this paper disincentivizes untruthfulness when reporting value of time. An approach grounded in mechanism design will help shed light on this analytically.  We have not seen this problem attacked in the lane-changing literature, but similar ideas have appeared recently in relation to bottleneck trading \citep{wang2018trading} and platooning of connected vehicles \citep{sun2019behaviorally}.  That the outcome of the games depends on communicated VOTs creates opportunities to \textit{game the system} and  it might be possible to engage in multiple transactions and exploit arbitrage opportunities.  Hence, appropriate pricing schemes would be a particularly interesting avenue to investigate.

\bibliography{refs}

\end{document}